\documentclass[12pt,leqno]{amsart}
\usepackage{amssymb}
\usepackage{amsmath,amssymb,color}
\numberwithin{equation}{section}
\newcommand{\DDD}{\mathcal{D}}

\newcommand{\ep}{\varepsilon}
\newcommand{\la}{\lambda}

\newcommand{\ppp}{\partial}
\newcommand{\SIGZZ}{\Sigma_{\gamma-\frac{\pi}{2}}}
\newcommand{\ABSz}{\vert z\vert}

\newcommand{\pppb}{\partial_t^{\beta}}

\newcommand{\rrrr}{\longrightarrow}
\newcommand{\ddda}{d_t^{\alpha}}
\newcommand{\AAA}{\mathcal{A}}
\newcommand{\ABAB}{(-A)^{\beta}}
\newcommand{\AEAE}{(-A)^{\ep}}

\newcommand{\sumk}{\sum_{k=0}^\infty}

\newcommand{\pppa}{\partial_t^{\alpha}}

\newcommand{\R}{\mathbb{R}}
\newcommand{\C}{\mathbb{C}}
\newcommand{\N}{\mathbb{N}}

\newcommand{\www}{\widetilde}

\newcommand{\ooo}{\overline}
\newcommand{\OOO}{\Omega}

\allowdisplaybreaks

\begin{document}
\title
[]
{
Operator approach for time-fractional evolution  
equations in Banach spaces}
%


\pagestyle{myheadings}

\author{
$^1$ Giuseppe Floridia, $^2$ Fikret G\"olgeleyen $^3$ Masahiro Yamamoto$^*$}
\thanks{
$^1$ 
Sapienza Universit\`a di Roma, 
Dipartimento di Scienze di Base e Applicate per l'Ingegneria
via Antonio Scarpa 16, 00161 Roma Italy
\\
$^2$ 
Department of Mathematics, Faculty of Science, 
Zonguldak B\"ulent Ecevit University, Zonguldak 67100 T\"urkiye
\\
$^3$ $^*$ Correspondence author\\
Graduate School of Mathematical Sciences, The University
of Tokyo, Komaba, Meguro, Tokyo 153-8914 Japan 
\\
Department of Mathematics, Faculty of Arts and Sciences, 
Zonguldak B\"ulent Ecevit University, Zonguldak 67100 Turkey\\
*Correspondence to [{\tt myama@ms.u-tokyo.ac.jp}]
}

\date{}

\begin{abstract}
Our first main purpose is to establish a framework for initial value problems 
for time-fractional evolution equation of order $\alpha \in (0,1)$
in Banach space $X$:
$$
\pppa (u(t)-a) = Au(t) + F(t), \quad 0<t<T.          \eqno{(*)}
$$
Here $u: (0,T) \rrrr X$ is an $X$-valued function defined in $(0,T)$,
and $a \in X$ is an initial value.
The operator $A$ satisfies a decay condition of resolvent which is 
the same as a generator of analytic semigroup.
Based on $X$-valued Laplace transforms, we establish a solution formula
yielding the well-posedness for (*).
In particular, we can directly treat 
a case $X=L^p(\OOO)$ over a bounded domain $\OOO$
and a uniform elliptic operator $A$.
Our theory is feasibly applicable to other topics 
such as regularity of solutions, inverse problems and control problems.
\end{abstract} 
\baselineskip 18pt

\maketitle
\section{Introduction and formulation}

Let $X$ be a Banach space over $\C$ and let $A: \DDD(A) \longrightarrow
X$ be a densely defined closed linear operator satisfying 
{\it Condition ($\AAA$)} stated later.  

For $u=u(t): [0,T] \rrrr X$, we introduce an initial value problem, called 
a forward problem, for a time-fractional evolution equation in $X$:
$$
\left\{ \begin{array}{rl}
& \ddda u(t) = Au(t) + F(t) \quad \mbox{for $t>0$ in $X$,}\\
& u(0) = a.
\end{array}\right.
                                          \eqno{(1.1)}
$$
Here we set $u'(t):= \frac{du}{dt}(t)$ and for $0<\alpha<1$, 
the Caputo derivative $\ddda u$ is defined by 
$$
\ddda u(t) = \frac{1}{\Gamma(1-\alpha)}\int^t_0 (t-s)^{-\alpha}
u'(s) ds,
$$
provided that the right-hand side can be defined (e.g., $u \in C^1([0,T];X)$).

Our main motivation is not only for studying forward problems such as
initial boundary value problem (1.1) as well as applications such as 
inverse problems for (1.1) in a Banach space $X$.
Needless to say, for applications such as inverse problems,   
we need some adequate framework and theory for the forward problem.
Thus we have to exploit
a convenient theory, not only with the scope limited to the forward 
problem but also to inverse problems, etc.  

There are many works on the well-posedness in the case where $X$ is 
a Hilbert space and $A$ admits an orthonormal basis composed of 
the eigenvectors of $A$ (for example, $A$ is self-adjoint with compact 
resolvent), and we can refer to Achache \cite{Ach1}, \cite{Ach2},
Jin \cite{J}, 
Kubica, Ryszewska and Yamamoto \cite{KRY}, Kubica and Yamamoto \cite{KY},
Zacher \cite{Za2009}.  
Here we do not intend any comprehensive references,
and see also the references therein.

Moreover, the forward problems themselves for (1.1) for the case of 
Banach space $X$ and more general systems
have been well-studied. 
We do not intend to provide a comprehensive list of works, and 
we refer, for example, to 
Bazhlekova \cite{Ba}, Cl\'ement and Londen \cite{ClL}, 
Cl\'ement, Londen and Simonett \cite{ClLS}, Gal and Warma \cite{GW}, 
Kim, Kim and Lim \cite{KKL}, Li and Li \cite{LiLi},
Pr\"uss and Simonett \cite{PS} 
Zacher \cite{Za2005} - \cite{Za2019}.  
Especially Da Prato and Iannelli \cite{DaP1}, 
Da Prato, Iannelli and Sinestrari \cite{DaP2} are pioneering articles
about forward problems for evolution integral equations including (1.1), and
Pr\"uss \cite{Pr} is a comprehensive source book for evolutionary integral
equations.  Lizama \cite{Liz1}, \cite{Liz2}, Zacher \cite{Za2019} are 
valuable surveys.

Here we should emphasize that one can choose other formulation 
of initial value problems
for time-fractional differential equations.
According to the references, the formulations and 
the definitions are not the same, so that 
we should be careful for the comparison of the results.
One conventional formulation of time-fractional 
evolution equations is based on a Volterra integral equation
and we can especially refer to Bazhlekova \cite{Ba}, Pr\"uss \cite{Pr}
and the references therein. As is explained in Section 6,  
the formulation by a Volterra equation is not convenient for 
application problems such as inverse problems of determining 
an initial value $a$ and/or a source term $F(t)$.
Moreover, there are very few works on inverse problems in Banach space 
$X$ and we seriously need a theory 
for the forward problem (1.1) which should be  
consistently applicable to inverse problems and various control 
problems.
Pr\"uss and Simonett \cite{PS} (pp.189-192)
treats time-fractional 
equations as an application of a general theory developed there.
However, the arguments are restricted to some range of orders 
$\alpha$, which is assumed to be specified large according to
the power $p$ of a state space $L^p$, and we do not know whether 
their approach is possible for $\alpha>0$ which is small to be out of
the range.

In this article, we will establish a theory of the forward problem 
for (1.1) in 
Banach space $X$ which should be directly applicable to 
inverse problems and control problems.
In the succceding article, based on the results in this
article, we will prove, for example,
uniqueness results for an inverse problem
of determining an initial value, which cannot be found 
in any existing works.
As by-product, we can prove the decay of solution as $t\to \infty$
in Theorem 2.1 (5) stated in Section 2.

In order to establish the theory of the forward problem, we modify 
a standard method for construction of analytic semigroups
via the Laplace transform for the case of the order $\alpha=1$.  
As monographs for $\alpha=1$, we can refer for example to 
Kato \cite{Ka}, Pazy \cite{Pa}, Tanabe \cite{Ta}, Yagi \cite{Ya}.
We emphasize that our methods for applications are very similar to 
the case $\alpha=1$ 
(e.g., Henry \cite{H}, Pazy \cite{Pa}, Yagi \cite{Ya}), and 
detailed studies should be future works.
%

Now, we formulate an initial boundary value problem for the time-fractional 
evolution equation.
To this end, 
we need to define fractional derivatives in 
the Lebesgue spaces and introduce a class of operators $A$ and function spaces.
Throughout this article, we assume that the indexes of the function spaces
satisfy 
$$
1 < p, q < \infty,
$$
if we do not specify.
Let $L^q(0,T;X)$ with $q\ge 1$ be the Lebesgue space of $X$-valued functions 
from $(0,T)$ and let
$$
\Vert v\Vert_{L^q(0,T;X)}:= \left( \int^T_0 \Vert v(s)\Vert_X^q ds \right)
^{\frac{1}{q}}.
$$
We will define $\pppa$ in $L^q(0,T;X)$ where $X$ is a Banach space $X$
(e.g., Yamamoto \cite{Y22}, \cite{Y25}).

For $\beta > 0$, we set 
$$
J^{\beta}v(t):= \frac{1}{\Gamma(\beta)}\int^t_0 (t-s)^{\beta-1}
v(s) ds, \quad v\in L^q(0,T;X).            \eqno{(1.2)}
$$
Then, we can easily prove that
$J^{\beta}: L^q(0,T;X) \longrightarrow J^{\beta}L^q(0,T;X)$ is injective and 
surjective.  We define a function space  
$$
W_{\beta,q}(0,T;X) := J^{\beta}L^q(0,T;X)                \eqno{(1.3)}
$$
with the norm $\Vert v\Vert_{W_{\beta,q}(0,T;X)}
:= \Vert (J^{\beta})^{-1}v \Vert_{L^q(0,T;X)}$.
We can verify that $W_{\beta,q}(0,T;X)$ is a Banach space.

We define a time-fractional derivative $\ppp_t^{\beta}$ 
in $W_{\beta,q}(0,T;X)$ by 
$$
\ppp_t^{\beta} := (J^{\beta})^{-1}, \quad
\DDD(\ppp_t^{\beta}) = J^{\beta}L^q(0,T;X).
                                                   \eqno{(1.4)}
$$
Henceforth we denote $J^{-\beta}:= (J^{\beta})^{-1}$ which means not only 
algebraically but topologically.
We understand that $\DDD(\cdot)$ is the domain of the operator 
under consideration.
Then 
\\
{\bf Lemma 1.1.}
\\
{\it
$\pppb : \DDD(\pppb) = W_{\beta,q}(0,T;X)\, \longrightarrow
\, L^q(0,T;X)$ is a closed operator.
}
\\
{\bf Proof.}
Let $u_n \in \DDD(\pppb)$ for $n\in \N$ and $u_n \longrightarrow u$ in 
$L^q(0,T;X)$, $\pppb u_n \longrightarrow w$ in $L^q(0,T;X)$ with some 
$u, w \in L^q(0,T;X)$.
Since $u_n\in \DDD(\pppb)$, we can find $w_n \in L^q(0,T;X)$ such that 
$u_n = J^{\beta}w_n$, that is, $w_n = \pppb u_n$,
Then $w_n \longrightarrow w$ and $J^{\beta} w_n \longrightarrow u$
in $L^q(0,T;X)$.
Since $J^{\beta}: L^q(0,T;X) \longrightarrow L^q(0,T;X)$ is bounded 
by the Young inequality on the convolution, we know that 
$J^{\beta} w_n \longrightarrow J^{\beta}w = u$ in $L^q(0,T;X)$.
This means that $w = \pppb u$ and $u \in \DDD(\pppb)$.
The proof of Lemma 1.1 is completed.
$\blacksquare$
\\

Henceforth, by $\rho(A)$ we denote the resolvent set of an operator 
$A$ and we write
$$
\Sigma_{\gamma}:= \{ z\in\C;\, \vert \mbox{arg}\, z\vert < \gamma,
\,\, z\ne 0\}
$$
with $0 < \gamma < \pi$.
 
For the operator $A$, we pose 
\\
{\bf Condition ($\AAA$):}
\\
(i) $A$ is a closed linear operator in $X$ such that $\DDD(A)$ is dense
in $X$.
\\
(ii) There exists a constant $\gamma' \in \left( \frac{\pi}{2}, \, \pi\right)$
such that $\Sigma_{\gamma'} \subset \rho(A)$.
\\
(iii) For arbitrarily fixed $\gamma \in (0, \gamma')$, there exists a constant 
$C = C_{\gamma}>0$ such that 
$$
\Vert (\la - A)^{-1}\Vert \le \frac{C_{\gamma}}{\vert \la\vert}
\quad \mbox{for all $\la \in \Sigma_{\gamma}$}.             \eqno{(1.5)}
$$
\\
(iv) $0 \in \rho(A)$.
\\

Then, the fractional power $(-A)^{\beta}$ with $0<\beta<1$ can be defined 
(e.g., \cite{Pa}, \cite{Ta}, \cite{Ya}).
In this article, we exculsively consider the case $0<\alpha<1$.

In terms of $\ppp_t^{\alpha}$ defined by (1.4), 
we formulate the initial value problem for
time-fractional evolution equation:
$$
\pppa (u(t)-a) = Au(t) + F(t) \quad \mbox{for $t>0$ in $X$.}
                                                            \eqno{(1.6)}
$$

The initial value problem (1.6) can describe an initial boundary 
value problem for a time-fractional diffusion equation.  More precisely,
let $\OOO \subset \R^d$ be a bounded domain with smooth boundary $\ppp\OOO$.
We formally define an elliptic operator $\mathcal{A}$ by
$$
\mathcal{A}v(x) := \sum_{k,\ell=1}^d \ppp_k(a_{k\ell}(x)\ppp_{\ell}v(x)) 
+ \sum_{k=1}^d b_k(x)\ppp_kv(x) + c(x)v(x),                        \eqno{(1.7)}
$$
where $a_{k \ell} = a_{\ell k} \in C^2(\ooo{\OOO})$, $b_k$, 
$c \in C(\ooo{\OOO})$, $c \le 0$ in $\OOO$, and we assume that 
there exists a constant $\kappa>0$ such that 
$$
\sum_{k,\ell=1}^d a_{k\ell}(x)\xi_k\xi_{\ell} 
\ge \kappa \sum_{k=1}^d \xi_k^2, \quad 
x\in \ooo{\OOO}, \, \xi_1, ..., \xi_d \in \R.    
$$
Henceforth $W^{2,p}(\OOO)$, etc. denotes usual Sobolev spaces
(e.g., Adams \cite{Ad}). 
In order to set up $\mathcal{A}$ within the framework of the
operator theory in Banach space $X$, to $\mathcal{A}$ we attach the domain 
$$
\DDD(A) = \{ u \in W^{2,p}(\OOO);\, 
u\vert_{\ppp\OOO} = 0\},                       \eqno{(1.8)}
$$
and by $A$ we define such an elliptic operator $\mathcal{A}$ with the domain 
$\DDD(A)$.  We can similarly consider other boundary conditions 
such as the homogeneous Neumann condition and the 
homogeneous Robin condition. 

Then, the condition $c\le 0$ implies that $0 \in \rho(A)$.
For $1<p<\infty$, it is known (e.g., Theorem 3.2 (p.213) in \cite{Pa})
that $A$ satisfies {\it Condition ($\AAA$)}.

The system (1.6) is the basic formulation of the forward problem
throughout the article.
Our main methodology is based on construction of solution mapping:
$(a,F) \mapsto u(t)$ for (1.6) by means of the theory for 
$X$-valued Laplace transforms (e.g., Arendt, Batty, Hieber and
Neubrander \cite{ABHN}).   We note that our solution mapping 
$(a,0) \mapsto u(t)$ for $0<\alpha<1$ corresponds to the semigroup 
$e^{tA}$ for the case
$\alpha=1$, although our solution mapping has no semigroup properties.
\\

This article is composed of the seven sections.
\begin{itemize}
\item
Section 2. Solution formula.\\
We state the essence of our forward theory as Theorem 2.1, which 
is systematically applied in later Sections 6 - 8.
\item
Sections 3 - 5. Proofs of the results in Section 2.
\item
Section 6. Comparison of the formulation by a Volterra integral 
equation. 
\item
Section 7. Conclusions.
\end{itemize}

The main purpose is the forward theory, and applications to inverse problems
and control problems are discussed in a succeeding article.

\section{Main results and solution formula}

We consider
$$
\pppa (u(t) - a) = Au(t) + F(t) \quad \mbox{for $t>0$ in $X$.}    \eqno{(2.1)}
$$
Henceforth we define the Laplace transform of $u: [0,\infty) \rrrr X$ by 
$$
(Lu)(\la):= \int^{\infty}_0 e^{-\la t} u(t) dt,
$$
provided that the right-hand side exists for $\la \in \C$.
We recall that $\gamma \in \left( \frac{\pi}{2}, \, \pi\right)$ is 
specified in {\it Condition ($\AAA$)}.
\\

We set
$$
G(z)a:= \frac{1}{2\pi i}\int_{\Gamma} 
e^{\la z}\la^{\alpha-1}(\la^{\alpha}-A)^{-1}a d\la \quad 
\mbox{for $a \in X$ and $z\in \Sigma_{\gamma-\frac{\pi}{2}}$}.
                               \eqno{(2.2)}
$$
Here we can choose a path $\Gamma\subset \rho(A)$ as follows:
$\Gamma$ starts at $+\infty e^{-i\gamma}$ to reach near 
the origin $0\in \C$ and then avoiding $0$, surrounds $0$, and again
goes to $+\infty e^{i\gamma}$.
In surrounding $0$, the path remains so that $\{ \vert \mbox{Re}\, z\vert;\, 
z\in \Gamma\}$ is sufficiently small which implies  $\Gamma \subset \rho(A)$.  
Since $\Gamma$ does not intersect $\{ z\le 0\}$, we see that $\la^{\alpha-1}
(\la^{\alpha} - A)^{-1}a$ is holomorphic in a neighborhood of 
$\Gamma$.  Hence, by means of Cauchy's integral theorem, we can verify that 
$G(z)a$ is invariant under transformations of $\Gamma$ satisfying the 
above conditions.  For example, as $\Gamma$ we can choose the following 
two choices which provide the same value $G(z)a$:
\\
(a) $z$-dependent $\Gamma:= \Gamma_1 \cup \Gamma_2 \cup \Gamma_3$.
$$
\Gamma_1 = \Gamma_1(z) := \left\{\rho e^{-i\gamma};\, \rho 
> \frac{1}{\vert z\vert} \right\},
\quad
\Gamma_2:= \left\{\frac{1}{\vert z\vert}e^{i\theta};\, 
-\gamma \le \theta \le \gamma
\right\}, 
$$
$$
\Gamma_3:= \left\{\rho e^{i\gamma};\, \rho > \frac{1}{\vert z\vert} \right\}.
                                                     \eqno{(2.3)}
$$
\\
(b) $z$-independent $\Gamma:= \Gamma^1 \cup \Gamma^2 \cup \Gamma^3$.
$$
\Gamma^1 = \{\rho e^{-i\gamma};\, \rho > \ep\},
\quad
\Gamma^2:= \{ \ep e^{i\theta};\, -\gamma \le \theta \le \gamma\},
\quad
\Gamma^3:= \{\rho e^{i\gamma};\, \rho > \ep\},
                                                     \eqno{(2.4)}
$$
where the constant $\ep>0$ is sufficiently small.

In this article, we mainly choose $\Gamma$ defined by (2.3).

Henceforth for $z \in \Sigma_{\gamma-\frac{\pi}{2}}
:= \left\{ z\in\C;\, \vert \mbox{arg}\, z\vert < \gamma - \frac{\pi}{2},\,
\, z\ne 0\right\}$ and $\beta > 0$, we define the Riemann-Liouville
fractional integral operator $J^{\beta}$ by
$$
J^{\beta}w(z):= \frac{1}{\Gamma(\beta)}\int^z_0 (z-s)^{\beta-1} w(s) ds
                                               \eqno{(2.5)}
$$
for $w \in L^1(\overrightarrow{0z})$, 
where the integral is done along the segment $\overrightarrow{0z} \subset
\C$ directed from $0$ to $z$.
  
The advantages of our main result Theorem 2.1 
for the forward problem are as follows:
\begin{itemize}
\item
it provides a solution formula yielding the well-posedness in various
solution classes such as weak solution, mild solution, strong solution,
classical solution.  Application methods of the solution formula are the same
as in the case $\alpha=1$ (see e.g., Theorem 6.1). 
\item
it is well incorporated with adequate uniqueness results of solutions
(Theorem 2.2).  
\item
it admits systematic applications to inverse problems, etc.
\end{itemize}

{\bf Theorem 2.1.}
\\
{\it
We assume Condition ($\mathcal{A}$).  
Let $\Gamma \subset \C$ is defined by (2.3).
Then, the infinite integral (2.2) converges in $X$ for each 
$z \in \SIGZZ$ and $a\in X$.
Moreover,  
$$
K(z)a:= \frac{d}{dz}J^{\alpha}G(z)a, \quad z\in \SIGZZ, \, a\in X
                                                   \eqno{(2.6)}
$$
exists in $X$.
\\
{\bf (1)}
For  
$$
1 < q < \infty,
$$
let  
$$
\mu
\left\{ \begin{array}{rl}
= &0 \quad  \mbox{if $\alpha q < 1$}, \\
> &1 - \frac{1}{q\alpha} \quad \mbox{if $\alpha q \ge 1$}.
\end{array}\right.
                                           \eqno{(2.7)}
$$
and let 
$$
a \in \DDD((-A)^{\mu}) \quad \mbox{and}\quad
F \in L^q(0,T; \DDD((-A)^{\ep})) \quad \mbox{with some $\ep>0$}.
$$
Then,
$$
\left\{ \begin{array}{rl}
& u(t) = G(t)a + \int^t_0 K(t-s)F(s) ds \quad \mbox{for $0<t<T$, 
exists and},  \\
& \mbox{satisfies $u\in L^q(0,T;\DDD(A))$, $u-a \in W_{\alpha,q}(0,T;X)$ 
and (2.1)}.
\end{array}\right.
                                               \eqno{(2.8)}
$$
Moreover, we can find a constant $C>0$ such that 
$$
\Vert u\Vert_{L^q(0,T;\DDD(A))} + \Vert u-a \Vert_{W_{\alpha,q}(0,T;X)}
\le C(\Vert a\Vert_{\DDD((-A)^{\mu})} 
+ \Vert F\Vert_{L^q(0,T;\DDD((-A)^{\ep}))})      \eqno{(2.9)}
$$
for all $a\in \DDD((-A)^{\mu})$ and $F \in L^q(0,T;\DDD((-A)^{\ep}))$.
\\
{\bf (2) (Laplace transform)} 
For each $a\in X$, the function $G(z)a$ in $z>0$ is holomorphic 
in $z\in \SIGZZ$ and
$$
L(Ga)(\la) = \la^{\alpha-1}(\la^{\alpha}-A)^{-1}a \quad 
\mbox{for $\la \in \Sigma_{\gamma-\frac{\pi}{2}}$ and $a \in X$}.
                                                   \eqno{(2.10)}
$$
\\
{\bf (3)}  For $0\le \beta \le 1$, there exists a constant $C=C(\beta)>0$ 
such that 
$$
\Vert \ABAB G(z)a\Vert \le C\vert z\vert^{-\alpha\beta}\Vert a\Vert \quad 
\mbox{for all $z \in \SIGZZ$ and $a \in X$}
                                                \eqno{(2.11)}
$$
and
$$
\Vert \ABAB K(z)a\Vert \le C\vert z\vert^{\alpha(1-\beta)-1}\Vert a\Vert \quad 
\mbox{for all $z\in \SIGZZ$ and $a \in X$.}
                                                \eqno{(2.12)}
$$
Here we emphasize that the constant $C>0$ is uniform for 
all $z\in \SIGZZ$. 
\\
{\bf (4)} We have $u\in C([0,T];X)$, that is,
$$
\lim_{t'\to t} \Vert G(t')a - G(t)a \Vert = 0
\quad \mbox{for all $0\le t \le T$ and $a\in X$}.
$$
\\
{\bf (5) (Decay)} We assume that $0 \in \rho(A)$.
For $a\in \DDD((-A)^{\mu})$, let $u(t)$ be the solution to 
$$
\pppa (u-a) = Au(t), \quad t>0.
$$
Then, there exists a constant $C>0$ such that  
$$
\Vert Au(t)\Vert 
\le Ct^{-\alpha}\Vert a\Vert \quad \mbox{for all $t>0$ and $a \in X$}.
$$
In particular, $\Vert u(t)\Vert \le Ct^{-\alpha}\Vert a\Vert$ 
for all $t>0$ and $a \in X$.
}

We say that $u$ is a strong solution to (2.1) if 
$u\in L^q(0,T;\DDD(A))$ satisfies (2.1) and
$u-a \in W_{\alpha,q}(0,T;X)$, as is shown in part (1) of the 
theorem.

In Theorem 2.1 (4), we interpret $\lim_{t' \to t} = \lim_{t'\downarrow t}$ for 
$t=0$ and $\lim_{t' \to t} = \lim_{t'\uparrow t}$ for 
$t=T$.

Theorem 2.1 generalizes the results for the case of $q=2$ and $X=L^2(\OOO)$ 
(e.g., Jin \cite{J},
Kubica, Ryszewska and Yamamoto \cite{KRY}, Sakamoto and 
Yamamoto \cite{SY}, for example), and the eigenfunction expansions of the 
solutions works in the case of Hilbert space $X$.
The representation formula (2.8) is the core achievement of 
Theorem 2.1 for applications, and for activating (2.8), we need 
estimates (2.11) and (2.12). The estimates (2.11) and (2.12) for Hilbert space 
$X$ and self-adjoint 
$A$, can be proved also by the eigenfunction expansions (e.g., 
Gorenflo, Luchko and Yamamoto \cite{GLY}). 
\\
{\bf Remark 2.1.}
\\
For a special case $\beta \in \N$, the proof of estimate (2.11) is 
easier for positive number $z$ (Theorem 2.2 (ii) (p.57) 
and Theorem 3.1 (p.73)) in Pr\"uss \cite{Pr}.
However, the case 
$\beta \in (0,1)$ is important for the applications as is seen 
in Section 6 for example.
\\
{\bf Remark 2.2.}
\\
Theorem 2.1 asserts that for the existence of solution 
$u \in L^q(0,T;\DDD(A))$ satisfying 
$u-a \in W_{\alpha,q}(0,T;X)$, we have to assume a stronger 
regularity $F \in L^q(0,T;\DDD((-A)^{\ep}))$ with some $\ep > 0$.
For $q\ne 2$ and a general Banach space, 
we do not know whether we can choose 
$\ep=0$, that is, $F \in L^q(0,T;X)$ is sufficient.
In particular, in the case of $a=0$, we call the maximal regularity 
if $Au$ and $\pppa u$ possess the same regularity as $F$.
Thus, if $F \in L^q(0,T;X)$, implies that  
$Au, \pppa u \in L^q(0,T;X)$ for at least one solution, then
we have the maximal regularity in the space $L^q(0,T;X)$.
The maximal regularity is known for the case of $q=2$ and a Hilbert space 
$X$ (e.g., Kubica, Ryszewska and Yamamoto \cite{KRY},
Sakamoto and Yamamoto \cite{SY}, Yamamoto \cite{Y22}, Zacher \cite{Za2009}). 
On the other hand, for the space $X=L^p(\OOO)$ with $p\ne 2$, similar 
results are known (e.g., Zacher \cite{Za2005}, \cite{Za2019}) under 
restrictive conditions which exclude some set of 
pairs $(\alpha,\, p) \in (0,1) \times (1, \,\infty)$.
We furthe refer to Corollaries 4.6 and 4.7 (pp.56-57) in
Bazhlekova \cite{Ba}, Zacher \cite{Za2006}.
In Section 7, we will discuss the maximal regularity in the case
of the H\"older continuous space in $t$.
\\
{\bf Remark 2.3.}
\\
The decay estimate of solution $u$ is proved directly as by-product 
within a general setting including $X=L^p(\OOO)$.  
As for a special case $X = L^2(\OOO)$, for example,
Kubica, Ryszewska and Yamamoto \cite{KRY}, 
Vergara and Zacher \cite{VZ}, prove the same decay as in Theorem 2.1 (5).
The case of symmetric $A$ is discussed in \cite{VZ}, 
including a time-fractional $p$-Laplacian equation and a
porus medium equation, while \cite{KRY} is concerned with
linear but not necessarily
symmetric $A$.  Both works are concerned with the cases where 
the coefficients depend on $x$ and $t$.  
See Kemppainen, Siljander, Vergara and Zacher \cite{KSVZ}
in the case of $x$-space $L^p(\R^d)$. 
\\

Theorem 2.1 does not assert the uniqueness of the solution, so that 
$u$ given by (2.8) is one solution to (2.1).
Now we state the uniqueness under several conditions.
For the uniqueness of solution to (2.1),
it suffices to prove that if $u\in L^q(0,T;\DDD(A)) \cap 
W_{\alpha,q}(0,T;X)$ satisfies $\pppa u = Au$ for $0<t<T$ where 
$T>0$ is finite or $T=\infty$, then 
$u=0$ for $0<t<T$.  Indeed, for 
$\pppa (u_1-a)(t) = Au_1(t) + F(t)$ and $\pppa (u_2-a)(t) 
= Au_2(t) + F(t)$ for $0<t<T$, we set $u:= u_1-u_2$ to 
have $\pppa u = Au$ for $0<t<T$.
\\
\vspace{0.2cm}
\\
{\bf Theorem 2.2.}
\\
{\it
(i) If $u\in W_{\alpha,q}(0,T;X) \cap L^q(0,T;\DDD(A))$ for all $T>0$ and
$\Vert u(t)\Vert_{X} = O(t^m)$ with some $m>0$ as $t\to\infty$
and $\pppa u(t) = Au(t)$ for $t>0$, then $u=0$ for $t>0$.
\\
(ii) Let $X$ be a Hilbert space with scalar product $(\cdot,\cdot)_X$ and 
$q\ge 2$,
and there exists a constant $C > 0$ such that 
$$
(Av, v)_X \le C\Vert v\Vert_X^2 \quad \mbox{for all $v \in \DDD(A)$}.
                                                      \eqno{(2.13)}
$$
If $u\in W_{\alpha,q}(0,T;X) \cap L^q(0,T;\DDD(A))$ satisfies 
$\pppa u = Au$ in $(0,T)$, then $u=0$ in $(0,T)$.
\\
(iii) 
Let $X = L^p(\OOO)$ with $1<p<\infty$ and let 
$1 < q < \infty$, and let $A$ be defined by (1.7) and (1.8) where
$a_{k\ell}, b_k, c \in C^{\infty}(\ooo{\OOO})$ and $c\le 0$ in $\OOO$.
If $u \in W_{\alpha,q}(0,T;X) \cap L^q(0,T;\DDD(A))$ satisfies 
$\pppa u = Au$ in $(0,T)$, then $u=0$ in $\OOO\times (0,T)$.
}

An elliptic opetator $A$ defined by (1.7) attached with 
the Dirichlet or the Neumann or the Robin boundary condition, 
satisfies (2.13).
\\

Our approach produces solution formula (2.8), which is 
convenient for applications also in non-Hilbert space $X$, 
as are shown in the succeeding sections.
We conclude this section with the existence of 
$u$ to (2.1) omitting the assumption $0 \in \rho(A)$.
\\
{\bf Corollary 2.1.}
\\
{\it
We assume that there exists a constant $C_0>0$ such that 
$$
\mbox{$A_0:= A-C_0$ with $\DDD(A_0) = \DDD(A)$ satisfies 
{\it Condition ($\AAA$)}}
$$
$$
\mbox{but not necessarily $0 \in \rho(A)$.}                     \eqno{(2.14)}
$$
Moreover let (2.7) hold.
Then, for $a \in \DDD((-A_0)^{\mu})$ and $F \in L^q(0,T; \DDD((-A_0)^{\ep}))$,
the problem (2.1) possesses a strong solution 
$u \in L^q(0,T;\DDD(A_0))$ satisfying \\
$u-a \in W_{\alpha,q}(0,T;X)$ and
$$
\Vert u\Vert_{L^q(0,T;\DDD(A_0))} + \Vert u-a\Vert_{W_{\alpha,q}(0,T;X)}
\le C(\Vert a\Vert_{\DDD((-A_0)^{\mu})} + \Vert F \Vert
_{L^q(0,T; \DDD((-A_0)^{\ep}))}).                \eqno{(2.15)}
$$
}
\\
{\bf Remark 2.4.}
Unlike the case $\alpha=1$, we cannot shift the zeroth-order term
by a simple transform $u \mapsto ue^{-C_0t}$.
Thus for the case $0<\alpha<1$, we have to resort to the iteration 
method.
\\

\section{Proof of Theorem 2.1: case of $a\ne 0$ and $F=0$}

{\bf First Step: Proof of Theorem 2.1 (2) and (5), and (2.11) .}
We will prove 
\\
{\bf Lemma 3.1.}
\\
{\it 
(i) $AG(z)a = G(z)Aa$ for all $z \in \Sigma_{\gamma-\frac{\pi}{2}}$ and 
$a \in \DDD(A)$.
\\
(ii) For $0\le \beta \le 1$, 
there exists a constant $C = C_{\beta} > 0$ such that 
$$
\Vert (-A)^{\beta}G(z)a\Vert \le C\vert z\vert^{-\alpha\beta}\Vert a\Vert
\quad \mbox{for all $z \in \Sigma_{\gamma-\frac{\pi}{2}}$ and $a \in X$.}
$$
(iii) $G(z)a$ is holomorphic in $z\in \Sigma_{\gamma-\frac{\pi}{2}}$ for
$a \in X$.
}
\\
In view of Theorem 2.6.1 (p.84) in Arendt, Batty, Hieber and 
Neubrander \cite{ABHN}, once Lemma 3.1 is proved, we can finish the proof of 
Theorem 2.1 (2) and (2.11).  Furthermore, Lemma 3.1 (ii) yields
Theorem 2.1 (5).
\\
{\bf Proof of Lemma 3.1.}
\\
(i) In terms of the closedness of the operator $A$, we can verify
\\
{\bf Lemma 3.2.}
\\
{\it
For an $X$-valued function $V: \Gamma \rrrr X$, we assume that 
$V(\la) \in \DDD(A)$ for each $\la \in \Gamma$ and
$AV(\cdot) \in C(\Gamma)$, $\Vert AV(\cdot)\Vert \in L^1(\Gamma)$.
Then 
$$
\int_{\Gamma} V(\la) d\la \in \DDD(A), \quad
A\left( \int_{\Gamma} V(\la) d\la\right)
= \int_{\Gamma} AV(\la) d\la.
$$
}
\\

We apply Lemma 3.2 to $V(\la) := \frac{1}{2\pi i}e^{\la z}
\la^{\alpha-1}(\la^{\alpha} - A)^{-1}a$.
We note that $\Gamma$ is defined by (2.3).  Since 
$$
AV(\la)
:= \frac{1}{2\pi i} e^{\la z}\la^{\alpha-1} (\la^{\alpha}-A)^{-1}Aa,
$$ 
by (1.5) with $\ep=0$, we obtain
$$
\Vert A^jV(\la)\Vert \le \frac{1}{2\pi}\sup_{\la\in \Gamma}
\vert e^{\mbox{Re}\,(\la z)}\vert \vert \la\vert^{\alpha-1}
\frac{\Vert A^ja\Vert}{\vert \la\vert^{\alpha}}
\le C \sup_{\la\in \Gamma}
\vert e^{\mbox{Re}\,(\la z)}\vert \frac{\Vert A^ja\Vert}{\vert \la\vert},
\quad \la \in \Sigma_{\gamma}
$$
for $j=0,1$.

We estimate Re $(\la z)$ for $z \in \Sigma_{\gamma-\frac{\pi}{2}}$.  Then,
$$
z = \vert z\vert e^{i\psi}, \quad 
 -\gamma+\frac{\pi}{2} < \psi < \gamma-\frac{\pi}{2}.     \eqno{(3.1)}
$$
On $\Gamma_1 \cup \Gamma_3$, we have 
$\la = \vert \la\vert e^{\pm i\gamma} 
= \vert \la\vert(\cos \gamma \pm i \sin \gamma)$ 
and so
$$
\mbox{Re}\, (\la z) = \mbox{Re}\, (\vert \la\vert
\vert z\vert e^{i\psi}e^{\pm i\gamma})
= \vert \la\vert \vert z\vert \cos (\gamma \pm \psi).
$$
For $-\gamma+\frac{\pi}{2} < \psi < \gamma-\frac{\pi}{2}$, we can directly 
verify that $\frac{\pi}{2} < \gamma \pm \psi < \frac{3}{2}\pi$.
Therefore, we can choose a constant $\delta_0>0$ such that 
$\cos (\gamma \pm \psi) \le -\delta_0$ for all 
$\vert \psi \vert < \gamma-\frac{\pi}{2}$.

For $\la \in \Gamma_2$, we see that 
$\la = \frac{1}{\vert z\vert}(\cos \theta + i \sin \theta)$ for 
$-\gamma < \theta < \gamma$, and so 
$\la z = (\cos \theta + i\sin\theta)(\cos \psi + i\sin\psi)$, 
which implies
$$
\mbox{Re}\, (\la z) =  \cos\theta \cos\psi - \sin \theta \sin \psi
= \cos (\theta - \psi).
$$
Thus
$$
\left\{ \begin{array}{rl}
& \mbox{Re}\, (\la z) \le -\vert \la\vert \vert z\vert \delta_0 \quad 
\mbox{for all $\la \in \Gamma_1 \cup \Gamma_3$ and 
$z \in \Sigma_{\gamma-\frac{\pi}{2}}$}, \\
& \mbox{Re}\, (\la z) = \cos (\theta - \psi) \quad 
\mbox{for all $\la \in \Gamma_2$ and 
$z \in \Sigma_{\gamma-\frac{\pi}{2}}$}.
\end{array}\right.
                                               \eqno{(3.2)}
$$
In particular, we have
$$
\left\{ \begin{array}{rl}
& \mbox{Re}\, (\la t) \le -\vert \la\vert t \delta_0 \quad 
\mbox{for all $\la \in \Gamma_1 \cup \Gamma_3$ and $t > 0$}, \\
& \mbox{Re}\, (\la t) 
= \cos \theta \quad \mbox{for all $\la \in \Gamma_2$ and $t>0$}.
\end{array}\right.
                             \eqno{(3.3)}
$$
Therefore, for $t>0$, we obtain
$$
\Vert A^jV(\la)\Vert \le
\left\{ \begin{array}{rl}
& Ce^{-\vert \la\vert t\delta_0}\vert z\vert \Vert A^ja\Vert,\quad
\la \in \Gamma_1 \cup\Gamma_3,  \\
& Ce \vert z\vert \Vert A^ja\Vert,\quad 
\la \in \Gamma_2
\end{array}\right.
$$
and for any fixed $t>0$, we have $\Vert A^jV(\cdot)\Vert 
\in L^1(\Gamma)$.  The rest assumptions of Lemma 3.2 are readily 
verified, so that the proof of Lemma 3.1 (i) is completed.
$\blacksquare$
\\
{\bf Proof of Lemma 3.1 (ii).}
\\
In terms of {\it Condition ($\AAA$)}, we can define a fractional power 
$(-A)^{\beta}$ of $-A$ (e.g., \cite{Pa}, \cite{Ta}).
Then, 
\\
{\bf Lemma 3.3.}
\\
{\it
Let $0 \le \beta \le 1$.  Then the operator $(-A)^{\beta}$ is closed in 
$X$, and there exists a constant $C_{\beta} > 0$ 
such that 
$$
\Vert (-A)^{\beta}(A-\la)^{-1}a\Vert
\le C_{\beta}\vert \la\vert^{\beta-1}\Vert a\Vert
$$
if $\la \in \Sigma_{\gamma}$.
}
\\
The proof is found e.g., Corollary (p.39) in \cite{Ta}.
By (2.2) and the closedness of the operator $(-A)^{\beta}$, we can justify 
$$
(-A)^{\beta}G(z)a = \frac{1}{2\pi i}\int_{\Gamma} 
e^{\la z}\la^{\alpha-1}(-A)^{\beta}(\la^{\alpha}-A)^{-1}a d\la, \quad 
z\in \Sigma_{\gamma-\frac{\pi}{2}},
$$
and so (3.3) yields
\begin{align*}
& \Vert (-A)^{\beta}G(z)a \Vert 
\le C \int_{\Gamma} 
\vert e^{\la z}\vert \vert \la\vert^{\alpha-1}
\Vert (-A)^{\beta}(\la^{\alpha}-A)^{-1}a\Vert \vert d \la\vert  \\
=& C\left( \int_{\Gamma_1} + \int_{\Gamma_2} + \int_{\Gamma_3}\right) 
\vert e^{\mbox{Re}\, (\la z)}\vert \vert \la\vert^{\alpha-1}
\Vert (-A)^{\beta}(\la^{\alpha}-A)^{-1}a\Vert \vert d\la\vert\\
\le& C\int_{\Gamma_1} e^{-\vert \la\vert\vert z\vert \delta_0}
\vert \la\vert^{\alpha-1} \Vert (-A)^{\beta}(\la^{\alpha}-A)^{-1}a\Vert
\vert d\la\vert 
+ C\int_{\Gamma_2} e^{\cos(\theta-\psi)}
\vert \la\vert^{\alpha-1} \Vert (-A)^{\beta}(\la^{\alpha}-A)^{-1}a\Vert
\vert d\la\vert \\
+& C\int_{\Gamma_3} e^{-\vert \la\vert \vert z\vert \delta_0}
\vert \la\vert^{\alpha-1} \Vert (-A)^{\beta}(\la^{\alpha}-A)^{-1}a\Vert
\vert d\la\vert \\
=: &I_1(z) + I_2(z) + I_3(z).
\end{align*}
{\bf Estimation of $I_1(z)$.}
We set $\la = \rho e^{-i\gamma}$ with $\frac{1}{\vert z\vert} < \rho < \infty$.
Then $d\la = e^{-i\gamma}d\rho$ and $\vert d\la\vert = d\rho$.
Using Lemma 3.3, we obtain
$$
I_1(z) 
= C\int^{\infty}_{\frac{1}{\vert z\vert}} e^{-\rho \vert z\vert\delta_0}
\rho^{\alpha-1}C\rho^{\alpha(\beta-1)} \Vert a\Vert \, d\rho
= C\int^{\infty}_{\frac{1}{\vert z\vert}} e^{-\rho\vert z\vert \delta_0}
\rho^{\alpha\beta-1} d\rho \Vert a\Vert.
$$
Setting $\eta = \rho \vert z\vert$, we have $d\eta = \vert z\vert d\rho$, 
and so 
$$
 \int^{\infty}_{\frac{1}{\vert z\vert}} e^{-\rho \vert z\vert \delta_0}
\rho^{\alpha\beta-1} d\rho
= \int^{\infty}_1 e^{-\delta_0 \eta}\left( \frac{\eta}{\vert z\vert}\right)
^{\alpha\beta-1} 
\frac{1}{\vert z\vert} d\eta 
= \vert z\vert^{-\alpha\beta}\int^{\infty}_1 e^{-\delta_0\eta} 
\eta^{\alpha\beta-1} d\eta
=: C\vert z\vert^{-\alpha\beta}.
$$
Consequently,
$$
\vert I_1(z)\vert \le C\vert z\vert^{-\alpha\beta}\Vert a\Vert.
$$
\\
{\bf Estimation of $I_2(z)$.}
We set $\la = \frac{1}{\vert z\vert}e^{i\theta}$ where $\theta: - \gamma 
\longrightarrow \gamma$.  Then $d\la = \frac{1}{\vert z\vert}
ie^{i\theta}d\theta$,
$\vert \la \vert = \frac{1}{\vert z\vert}$ and $\vert d\la\vert 
= \frac{1}{\vert z\vert}d\theta$.
Hence,
\begin{align*}
& I_2(z) = C\int^{\gamma}_{-\gamma} 
e \left( \frac{1}{\vert z\vert}\right)^{\alpha-1}
C \left( \frac{1}{\vert z\vert}\right)^{\alpha(\beta-1)}
\frac{1}{\vert z\vert} d\theta \Vert a\Vert \\
=& \frac{Ce}{\vert z\vert^{\alpha\beta}} \int^{\gamma}_{-\gamma}
 d\theta \Vert a\Vert
\le C\vert z\vert^{-\alpha\beta} \Vert a\Vert
= \frac{2C\gamma e}{\vert z\vert^{\alpha\beta}}\Vert a\Vert.
\end{align*}
For $I_3(z)$, we can estimate similarly.
Thus, the proof of Lemma 3.1 is complete.

Part (iii) follows from Theorem 2.6.1 (p.84) in \cite{ABHN}, which asserts the 
equivalence between the holomorphy of $G(z)a$ and the Laplace transform
$L(G(z)a)(\la)$.
Thus the proof of Lemma 3.1 is complete.
$\blacksquare$
\\
{\bf Second Step: Laplace transform of the fractional derivative.}
\\
In this step, we prove
\\
{\bf Lemma 3.4.}
\\
{\it
Let $u \in W_{\alpha,q}(0,T;X)$ for any $T>0$ and $\Vert u(\cdot,t)\Vert
= O(t^m)$ as $t \to \infty$ with some $m\ge 0$.
Then, $(L\pppa u)(\la)$ exists for $\la > 0$ and
$$
(L\pppa u)(\la) = \la^{\alpha}(Lu)(\la) \quad \mbox{for $\la > 0$}.
                            \eqno{(3.4)}
$$
}
{\bf Proof of Lemma 3.4.}
\\
For arbitrary $T>0$, since $u\in W_{\alpha,q}(0,T;X)$ by the definition
(1.4), we can find
$w_T \in L^q(0,T;X)$ such that $u = J^{\alpha}w_T$ and 
$\pppa u = w_T$ in $(0,T)$.  Then, 
$$
J^{1-\alpha}u = J^{1-\alpha}J^{\alpha}w_T
= J^{1-\alpha+\alpha}w_T = J^1w_T = \int^t_0 w_T(\cdot,s) ds
$$
for $0<t<T$.  Hence, $\frac{d}{dt}(J^{1-\alpha}u) \in L^q(0,T;X)$,
that is, $\pppa u(t) = \frac{d}{dt}(J^{1-\alpha}u)(t)$ for
$0<t<T$.
Moreover, the H\"older inequality yields 
$$
\Vert J^{1-\alpha}u(t)\Vert \le \int^t_0 \Vert w_T(s)\Vert ds
\le \left( \int^t_0 ds\right)^{\frac{q-1}{q}}\left( \int^t_0
\Vert w_T(s)\Vert^q ds\right)^{\frac{1}{q}} 
\, \longrightarrow \, 0
$$
as $t \to 0$ by $q>1$.
Therefore, the integration by parts implies
\begin{align*}
& \Gamma(1-\alpha)\int^T_0 \pppa u(t)e^{-\la t} dt
= \Gamma(1-\alpha)\int^T_0 \frac{d}{dt}(J^{1-\alpha}u)(t)e^{-\la t} dt\\
= &\int^T_0 \frac{d}{dt}\left( 
\int^t_0 (t-s)^{-\alpha}u(s) ds\right) e^{-\la t} dt\\
=& \left[ \left(\int^t_0 (t-s)^{-\alpha}u(s) ds\right) e^{-\la t}
\right]^{t=T}_{t=0}
+ \int^T_0 \left( \int^t_0 (t-s)^{-\alpha}u(s) ds\right)
\la e^{-\la t} dt.
\end{align*}
Then
$$
\lim_{T\to \infty} \left(\int^T_0 (T-s)^{-\alpha}u(s) ds\right) e^{-\la T} = 0.
$$
Indeed, 
\begin{align*}
& \left\Vert \int^T_0 (T-s)^{-\alpha}u(s) ds \right\Vert e^{-\la T} 
\le \left( \int^T_0 (T-s)^{-\alpha}\Vert u(s)\Vert ds \right)e^{-\la T} \\
\le & C\left( \int^T_0 (T-s)^{-\alpha}T^m ds \right)e^{-\la T} 
= CT^me^{-\la T}\frac{T^{1-\alpha}}{1-\alpha} \longrightarrow 0
\end{align*}
as $T \to \infty$.

Moreover, 
$$
 \la\int^T_0 \left(\int^t_0 (t-s)^{-\alpha}u(s) ds\right) e^{-\la t} dt
= \la\int^T_0 \left(\int^T_s (t-s)^{-\alpha}e^{-\la t} dt\right) u(s) ds.
$$
Here setting $\eta:= t-s$, that is, $t = \eta + s$, we have
$$
\int^T_s (t-s)^{-\alpha}e^{-\la t} dt
= e^{-\la s}\int^{T-s}_0 \eta^{-\alpha}e^{-\la\eta} d\eta
\longrightarrow e^{-\la s}\frac{\Gamma(1-\alpha)}{\la^{1-\alpha}}
$$
as $T \to \infty$.
Hence,
$$
\la\int^T_0 \left(\int^t_0 (t-s)^{-\alpha}u(s) ds\right) e^{-\la t} dt
\, \longrightarrow \, \la^{\alpha}\Gamma(1-\alpha)(Lu)(\la).
$$
Thus the proof of Lemma 3.4 is complete.
$\blacksquare$
\\
{\bf Third Step: Proof of Theorem 2.1 (1) in the case of $F=0$.}
\\
This step is devoted to the proof of 
\\
{\bf Proposition 3.1.}
\\
{\it
We arbitrarily fix $T>0$ and a constant $\mu$ satisfying (2.7).
Then,
$$
Ga - a \in W_{\alpha,q}(0,T;X) \cap L^{\infty}(0,\infty;X) \quad
\mbox{for all $a \in \DDD((-A)^{\mu})$,}
                                                           \eqno{(3.5)}
$$
and
$$
\pppa (G(t)a-a) = AG(t)a \quad \mbox{for $a\in \DDD((-A)^{\mu})$ and
all $0 < t < T$}.
                                         \eqno{(3.6)}
$$
}
We note that (3.5) means $\pppa (Ga - a) \in L^q(0,T;X)$.
\\
{\bf Proof of Proposition 3.1.}
\\
We prove
$$
A(LGa)(\la) = L(AGa)(\la) \quad \mbox{for all $\la > 0$ and
$a\in \DDD(A)$.}                    \eqno{(3.7)}
$$
{\bf Verification of (3.7).}
\\
We arbitrarily fix $\la > 0$.  In terms of $a\in \DDD(A)$, Lemma 3.1 (i)
yields $G(t)a \in \DDD(A)$ and $AG(t)a = G(t)Aa$ for all $t>0$.
Moreover $A(G)a\in C([0,\infty);X) \cap L^{\infty}(0,\infty;X)$ and
$e^{-\la t}G(t)a \in \DDD(A)$ for each $t>0$.
Then $e^{-\la t}AG(t)a \in C_t([0,\infty);X) \cap L_t^1(0,\infty;X)$ for 
arbitrarily fixed $\la > 0$.  
Here $C_t([0,\infty);X)$ and $L_t^1(0,\infty;X)$ 
mean the corresponding function
spaces with the independent variable $t$.

A similar application of Lemma 3.2 yields  
\begin{align*}
& AL(Ga)(\la) = A\int^{\infty}_0 e^{-\la t}G(t)a dt
= \int^{\infty}_0 A(e^{-\la t}G(t)a) dt
= \int^{\infty}_0 e^{-\la t}(AG(t)a) dt\\
=& L(AGa)(\la) \quad \mbox{for all $\la > 0$.}
\end{align*}
Thus the verification of (3.7) is complete.
$\blacksquare$

{\it Condition ($\AAA$)} yields 
$$
\sup_{\la \in \Sigma_{\gamma}} \Vert \la \la^{\alpha-1}(\la^{\alpha} - A)^{-1}a
\Vert < \infty,
$$
and so we can apply Theorem 2.6.1 in \cite{ABHN} to conclude that 
$$
G(z)a = \frac{1}{2\pi i}\int_{\Gamma} e^{\la z} \la^{\alpha-1}
(\la^{\alpha} - A)^{-1} a d\la
$$
is holomorphic in $\Sigma_{\gamma-\frac{\pi}{2}}$ and
$$
(LGa)(\la) = \la^{\alpha-1}(\la^{\alpha}-A)^{-1}a \quad
\mbox{for all $\la > 0$}.        \eqno{(3.8)}
$$

By (3.8), we have
$$
 \la^{\alpha}L(Ga)(\la) - A(LGa)(\la)
= (\la^{\alpha}-A)L(Ga)(\la) = \la^{\alpha-1}a
$$
for $\la > 0$.  Therefore, (3.7) yields
$$
\la^{\alpha}L(Ga)(\la) - L(AGa)(\la) = \la^{\alpha-1}a
$$
for $\la > 0$.  Since $(La)(\la) = \la^{-1}a$ by direct calculations,
we obtain
$$
\la^{\alpha}L(Ga-a)(\la) = L(AGa)(\la), \quad \la > 0.
                                   \eqno{(3.9)}
$$
Setting $w(t) := AG(t)a$, we rewrite (3.9) as
$$
(Lw)(\la) = \la^{\alpha}L(Ga-a)(\la), \quad \la > 0.
                                       \eqno{(3.10)}
$$
In view of Lemma 3.1 (i) and (ii), we have
$$
w(t) = G(t)Aa \in L^{\infty}(0,\infty;X) \subset L^q(0,T;X) \quad
\mbox{for any $T>0$.}
$$
Hence, setting 
$$
v:= J^{\alpha}w \in W_{\alpha,q}(0,T;X) \quad \mbox{for any $T>0$},
                                  \eqno{(3.11)}
$$
we see that $w = \pppa v$ and
\begin{align*}
& \Vert v(t)\Vert_X = \left\Vert \frac{1}{\Gamma(\alpha)}
\int^t_0 (t-s)^{\alpha-1} w(s) ds\right\Vert \\
\le& C\int^t_0 (t-s)^{\alpha-1}\Vert w(s)\Vert ds 
= C\Vert w\Vert_{L^{\infty}(0,\infty;X)}t^{\alpha} \quad
\mbox{for all $t>0$.}
\end{align*}
Therefore, Lemma 3.4  yields
$$
(Lw)(\la) = (L(\pppa v))(\la) = \la^{\alpha}Lv(\la), \quad \la > 0.
                                                    \eqno{(3.12)}
$$
In terms of (3.10), we have
$$
\la^{\alpha}(Lv)(\la) = \la^{\alpha}L(Ga - a)(\la), \quad \la > 0,
$$
that is,
$$
(Lv)(\la) = L(Ga - a)(\la), \quad \la > 0.
$$
The injectivity of the Laplace transform yields
$$
v(t) = G(t)a - a, \quad t>0.            \eqno{(3.13)}
$$
Since $v \in W_{\alpha,q}(0,T;X)$ for all $T>0$, by (3.11) 
we have proved that $Ga - a \in W_{\alpha,q}(0,T;X)$ in (3.5).

Finally (3.12) and (3.13) yield
\begin{align*}
& (Lw)(\la) = (L\pppa v)(\la)
= (L(\pppa (Ga - a)))(\la)
= \la^{\alpha}(Lv)(\la)\\
=& \la^{\alpha}(L(Ga-a))(\la), \quad \la > 0.
\end{align*}
Hence, (3.9) implies
$$
L(\pppa (Ga-a))(\la) = L(AGa)(\la), \quad \la>0.
$$
Again the injectivity of the Laplace transform yields
$AG(t)a = \pppa (G(t)a - a)$ for $t>0$.
Thus the proof of Proposition 3.1 is completed for $a \in \DDD(A)$.
$\blacksquare$

Finally, we will prove (3.5) and (3.6) for $a \in \DDD((-A)^{\mu})$, where
$\mu$ satisfies (2.7).
Since $\DDD(A)$ is dense in $X$, for each $a \in \DDD((-A)^{\mu})$, setting
$b:= (-A)^{\mu}a\in X$, we can find a sequence $b_n\in \DDD(A)$, $n\in \N$ 
such that $b_n \rrrr b$ in $X$ as $n \to \infty$.
Therefore, $(-A)^{-\mu}b_n \rrrr (-A)^{-\mu}b = a$ in $\DDD((-A)^{\mu})$
as $n\to \infty$.
By $b_n \in \DDD(A)$, we have 
$A((-A)^{-\mu}b_n) = (-A)^{-\mu}(Ab_n) \in X$, that is,
$a_n:= (-A)^{-\mu}b_n \in \DDD(A)$.  Hence, for each $a \in \DDD((-A)^{\mu})$,
we can find a sequence $a_n \in \DDD(A)$, $n\in \N$ such that 
$(-A)^{\mu}a_n \rrrr (-A)^{\mu}a$ in $X$ as $n\to \infty$, that is,
$$
a_n \rrrr a \quad \mbox{in $\DDD((-A)^{\mu})$ as $n\to \infty$}.
                                         \eqno{(3.14)}
$$

As is already proved, since $a_n \in \DDD(A)$, 
we see that $G(t)a_n - a_n$ satisfies 
(3.5) and (3.6).  By (3.14) and Lemma 3.1 (ii), we can obtain
$$
\Vert A(G(t)(a_n-a))\Vert = \Vert (-A)^{1-\mu}G(t)(-A)^{\mu}(a_n-a)\Vert
\le Ct^{-\alpha(1-\mu)}\Vert (-A)^{\mu}(a_n-a)\Vert.
$$
Consequently,
$$
\Vert AG(a_n-a)\Vert_{L^q(0,T;X)}
\le C\left( \int^T_0 t^{-\alpha q(1-\mu)} dt \right)^{\frac{1}{q}}
\Vert (-A)^{\mu}(a_n-a)\Vert
\le C_1\Vert (-A)^{\mu}(a_n-a)\Vert.
$$
Here we used
$$
\int^T_0 t^{-\alpha q(1-\mu)} dt = \frac{T^{1-\alpha q(1-\mu)}}
{1-\alpha q(1-\mu)} < \infty
$$
by (2.7).  Therefore, $AGa_n \rrrr AGa$ in $L^q(0,T;X)$ as 
$n\to \infty$.

Since $\pppa (Ga_n-a_n) = SAGa_n$ by (3.6), we see that 
$$
\pppa (Ga_n-a_n) \rrrr AGa, \quad
Ga_n - a_n \in \DDD(\pppa).
                             \eqno{(3.15)}
$$
Moreover, 
$$
\Vert (Ga_n - a_n) - (Ga-a)\Vert_{L^q(0,T;X)}
\le \Vert G(a_n - a)\Vert_{L^q(0,T;X)}
+ \Vert a_n-a\Vert_{L^q(0,T;X)} \, \rrrr \, 0
                                          \eqno{(3.16)}
$$
as $n\to \infty$ by means of Lemma 3.1 (ii) with $\beta = 0$.
Since $\pppa : L^q(0,T;X) \rrrr L^q(0,T;X)$ is a closed 
operator, the limits (3.15) and (3.16) yield 
$$
Ga - a \in \DDD(\pppa), \quad \mbox{and}\quad
\pppa (G(t)a - a) = AG(t)a.
$$
By means of Lemma 3.1 (ii), we can readily derive that 
$Ga - a \in L^q(0,T;X)$ for all $a\in \DDD((-A)^{\mu})$.
Thus the proof of Proposition 3.1 is completed for arbitrary $a \in 
\DDD((-A)^{\mu})$.
$\blacksquare$
\\
{\bf Fourth Step: Proof of Theorem 2.1 (4).}
\\
We consider only the case $0<t<T$ and the arguments for
the cases $t=0$ and $t=T$ are the same. 
First we will prove the proposition for $a \in \DDD(A)$.

By (2.2), we have 
$$
G(t')a - G(t)a
= \frac{1}{2\pi i}\int_{\Gamma} (e^{\la t'} - e^{\la t})
\la^{\alpha-1}(\la^{\alpha} - A)^{-1}a d\la.
$$
We directly see
$$
\lim_{t'\to t}(e^{\la t'} - e^{\la t}) \la^{\alpha-1}(\la^{\alpha} - A)^{-1}a 
= 0 \quad \mbox{for all fixed $\la \in \Gamma$}.
$$
Moreover, by (3.3) we obtain
$$
\vert e^{\la t'}\vert, \,\, \vert e^{\la t}\vert
\le
\left\{ \begin{array}{rl}
& e^{-\vert \la\vert t\delta} + e^{-\vert \la\vert t'\delta}, 
\quad \la \in \Gamma_1 \cup \Gamma_3, \\
& e, \quad \la \in \Gamma_2.
\end{array}\right.
$$
Hence Condition ($\mathcal{A}$) yields 
$$
\Vert (e^{\la t'} - e^{\la t}) \la^{\alpha-1}(\la^{\alpha} - A)^{-1}a \Vert
\le
\left\{ \begin{array}{rl}
& Ce^{-\vert \la\vert\min\{t, t'\}\delta}\frac{1}{\vert \la\vert}, \quad
\la \in \Gamma_1 \cup \Gamma_3, \\
& \frac{C}{\vert \la\vert}, \quad \la \in \Gamma_2.
\end{array}\right.
                                            \eqno{(3.17)}
$$
By recalling the definition (2.3) of $\Gamma$, the right-hand side of (3.17) 
is in $L^1(\Gamma)$.  Therefore, the Lebesgue convergence theorem completes 
the proof of Theorem 2.1 (4).
Thus we finish the proofs of the parts of Theorem 2.1 which are not 
concerned with $K(t)$.
$\blacksquare$

\section{Proof of Theorem 2.1: case of $a=0$ and $F\ne 0$}

In this section, we consider the case where the initial value $a$ is zero and 
the non-homogeneous term $F$ is not zero:
$$
\pppa u(t) = Au(t) + F(t),   \quad t>0.              \eqno{(4.1)}
$$
\\
{\bf First Step: Estimation of $\ABAB \frac{dG}{dz}(z)a$.}
We recall 
$$
G(z)a = \frac{1}{2\pi i}\int_{\Gamma} e^{\la z}\la^{\alpha-1}
(\la^{\alpha} - A)^{-1}a d\la,
$$
where the path $\Gamma \subset \C$ is defined by (2.3).

For any fixed $z \in \SIGZZ$, we see that 
$$
\frac{\ppp}{\ppp z}(e^{\la z}\la^{\alpha-1}(-A)^{\beta}
(\la^{\alpha}-A)^{-1}a) \in L^{\infty}(\Gamma) \cap L^1(\Gamma)
\quad \mbox{as a function in $\la$.}
$$
In terms of the Lebesgue convergence theorem, 
we can justify the exchange of $\frac{\ppp}{\ppp t}$ and 
$\frac{1}{2\pi i} \int_{\Gamma} \cdots d\la$, so that  
$$
\frac{\ppp}{\ppp z}(G(z)a)
= \frac{1}{2\pi i}\int_{\Gamma} \frac{\ppp}{\ppp z}
(e^{\la z}\la^{\alpha-1}(\la^{\alpha} - A)^{-1}a) d\la.
                                                        \eqno{(4.2)}
$$
Let 
$$
0 \le \beta \le 1.
$$
Then, 
\begin{align*}
& (-A)^{\beta}\frac{dG}{dz}(z)a
= \frac{1}{2\pi i}\int_{\Gamma} e^{\la z}\la^{\alpha}(-A)^{\beta}
(\la^{\alpha} - A)^{-1}a d\la\\
=& \frac{1}{2\pi i}\left( \int_{\Gamma_1} + \int_{\Gamma_2}
+ \int_{\Gamma_3}\right)  e^{\la z}\la^{\alpha}(-A)^{\beta}
(\la^{\alpha} - A)^{-1}a d\la
=: K_1(z) + K_2(z) + K_3(z).
\end{align*}
\\
{\bf Estimation of $K_1(z)$.}
We set $\la = \rho e^{-i\gamma}$ with $\rho > \frac{1}{\vert z\vert}$.
In view of Lemma 3.3, using (3.3), we have
\begin{align*}
& \Vert K_1(z)\Vert = \frac{1}{2\pi}\int_{\Gamma_1}
 e^{-\rho \ABSz \delta_0} \vert \la\vert^{\alpha}
\Vert (-A)^{\beta}(\la^{\alpha} - A)^{-1}a\Vert \vert d\la\vert\\
\le& C\int^{\infty}_{\frac{1}{\ABSz}} e^{-\rho \ABSz\delta_0} \rho^{\alpha}
\rho^{\alpha(\beta-1)} d\rho \Vert a\Vert
\le C\int^{\infty}_{\frac{1}{\ABSz}} e^{-\rho \ABSz\delta_0} 
\rho^{\alpha\beta} d\rho \Vert a\Vert.
\end{align*}
Changing the variables $\rho \mapsto \eta$ in the integral 
$\rho = \frac{\eta}{\ABSz}$,
we calculate
$$
  \int^{\infty}_{\frac{1}{\ABSz}} e^{-\rho \ABSz\delta_0} \rho^{\alpha\beta} 
d\rho
= \ABSz^{-\alpha\beta-1}\int^{\infty}_1 e^{-\delta_0\eta} 
\eta^{\alpha\beta} d\eta
=: C_{\beta}\ABSz^{-\alpha\beta - 1},
$$
which implies 
$$
\Vert K_1(z)\Vert \le C_{\beta}\frac{1}{\ABSz^{\alpha\beta+1}}\Vert a\Vert.
$$
\\
{\bf Estimation of $K_2(z)$.}
We set $\la = \frac{1}{\ABSz}e^{i\theta}$, where $\theta$ varies from 
$-\gamma$ to $\gamma$.  Then, $d\la = \frac{1}{\ABSz}ie^{i\theta}d\theta$ and 
$\vert \la\vert = \frac{1}{\ABSz}$, $\vert d\la\vert = \frac{1}{\ABSz}d\theta$.
Hence,
\begin{align*}
& \Vert K_2(z)\Vert \le C\int^{\gamma}_{-\gamma}
e^{\cos \theta} \left( \frac{1}{\ABSz}\right)^{\alpha}
\vert \la\vert^{\alpha(\beta-1)}\Vert a\Vert \frac{1}{\ABSz} d\theta\\
\le & C\int^{\gamma}_{-\gamma}
\frac{1}{\ABSz^{\alpha\beta+1}} \Vert a\Vert d\theta
\le C\ABSz^{-\alpha\beta-1} \Vert a\Vert.
\end{align*}

Similarly to $K_1(z)$, we can estimate $K_3(z)$ to obtain
\\
{\bf Lemma 4.1}
\\
{\it 
For $0\le \beta \le 1$, there exists a constant $C_{\beta} > 0$ such that 
$$
\left\Vert (-A)^{\beta}\frac{d}{dz}G(z)a\right\Vert \le C\ABSz^{-\alpha\beta-1}
\Vert a\Vert \quad \mbox{for all $z \in \SIGZZ$ and $a \in X$.}
$$
}
\\
{\bf Second Step: Estimation of $(-A)^{\beta}J^{\tau}\frac{dG}{dz}(z)a$}. 
\\
We will prove
\\
{\bf Lemma 4.2.}
\\
{\it
Let $0 \le \beta \le 1$ and $0 < \tau \le 1$.  Then, we have
$$
\frac{d}{dz}J^{\tau}Ga\in L^1(0,T;X), \quad  
\left\Vert (-A)^{\beta} \frac{d}{dz}J^{\tau}G(z)a\right\Vert 
\le C\ABSz^{(\tau-\alpha\beta)-1}\Vert a\Vert
$$
for all $z\in \SIGZZ$ and $a\in X$.
\\
In particular, choosing $\tau = \alpha$, we have
$$
\Vert (-A)^{\beta} K(z)a\Vert 
\le C\ABSz^{\alpha(1-\beta)-1}\Vert a\Vert \quad \mbox{for all 
$a \in X$ and $z\in \SIGZZ$}.
$$
}
Thus Lemma 4.2 completes the proof of (2.12).
\\
{\bf Proof of Lemma 4.2.}
\\
For $0 < \tau \le 1$, we note  
$$
\frac{d}{dz}(J^{\tau}v)(z) 
= \frac{1}{\Gamma(\tau)}\frac{d}{dz}
\int^z_0 (z-s)^{\tau-1}v(s) ds \quad \mbox{for $v\in L^1(0,T;X)$}.
$$

Henceforth, for $\alpha_1, \alpha_2 > 0$, we define the Mittag-Leffler 
function by
$$
E_{\alpha_1,\alpha_2}(z):= \sum_{k=0}^\infty \frac{z^k}
{\Gamma(\alpha_1 k + \alpha_2)}, \quad z\in \C
$$
and we know that the radius of convergence of the power series is $\infty$
and $E_{\alpha_1,\alpha_2}(z)$ is an entire function in $z \in \C$
(e.g., Podlubny \cite{Po}).

First we prove
\\
{\bf Lemma 4.3.}
\\
{\it 
Let $0 < \tau \le 1$ and $\la \in \C$.  Then,
$$
\frac{d}{dz}(J^{\tau}e^{\la z})(z)
= z^{\tau-1}E_{1,\tau}(\la z) \quad \mbox{for Re}\, z >0.
$$
}
\\
{\bf Proof of Lemma 4.3.}
\\
We have
$$
 (J^{\tau}e^{\la z})(z)
= \frac{1}{\Gamma(\tau)}\int^z_0 (z-s)^{\tau-1}
\sumk \frac{\la^ks^k}{k!} ds
= \frac{1}{\Gamma(\tau)}\sumk \frac{\la^k}{k!}
\int^z_0 (z-s)^{\tau-1}s^k ds
$$
$$
= \frac{1}{\Gamma(\tau)}\sumk \frac{\la^k}{k!}
\frac{\Gamma(\tau)\Gamma(k+1)}{\Gamma(\tau+k+1)}
z^{\tau+k}
= \sumk \frac{1}{\Gamma(\tau+k+1)} \la^k z^{k+\tau}
= z^{\tau}E_{1,\tau+1}(\la z),           \eqno{(4.3)}
$$
so that 
\begin{align*}
& \frac{d}{dz}(J^{\tau}e^{\la z})(z)
= \sumk \frac{\la^k(k+\tau)z^{k+\tau-1}}
{\Gamma(\tau+k+1)}
= \sumk \frac{(\tau+k)\la^kz^k}{(\tau+k)\Gamma(\tau+k)}
z^{\tau-1}\\
= & z^{\tau-1}\sum_{k=0}^{\infty} \frac{(\la z)^k}{\Gamma(\tau+k)}
= z^{\tau-1}E_{1,\tau}(\la z).
\end{align*}
Thus the proof of Lemma 4.3 is complete.
$\blacksquare$

Next, we show two lemmata.
\\
{\bf Lemma 4.4.}
\\
{\it
We choose $\sigma \in \left( \frac{\pi}{2},\, \pi\right)$ and 
$\tau \in (0, \,1]$ arbitrarily.  Then
$$
\vert E_{1,\tau}(\la t)\vert \le \frac{C}{1+\vert \la t\vert}
\quad \mbox{for all $t\ge 0$ and $\sigma \le \vert \mbox{arg}\, \la\vert 
\le \pi$}.
$$
}
The proof is found as Theorem 1.6 (p.35) in Podlubny \cite{Po}.

Recalling that $\frac{\pi}{2} < \gamma < \pi$, we can 
choose $\gamma_0$ such that $\frac{\pi}{2} < \gamma_0 < \gamma < \pi$.
Therefore, $\Gamma \subset \{ z\in \C;\, 
\gamma_0 \le \vert \mbox{arg}\, z\vert \le \pi\}$ and
$$
\vert E_{1,\tau}(\la t) \vert \le \frac{C}{1+\vert \la t\vert}
\quad \mbox{for $t\ge 0$ and $\la \in \Gamma$.}
                                                    \eqno{(4.4)}
$$

Now we show
\\
{\bf Lemma 4.5.}
\\
{\it
Let $0 \le \beta \le 1$.  Then there exists a constant 
$C=C_\beta > 0$ such that 
$$
\Vert J^{\beta}G(z)a\Vert \le C\ABSz^{\beta}\Vert a\Vert \quad
\mbox{for all $z \in \SIGZZ$ and $a\in X$}.
$$
In particular, 
$J^{\beta}Ga \in L^{\infty}(0,T;X)$.
}
\\
{\bf Proof of Lemma 4.5.}
\\
We have already proved (2.11), and so $\Vert G(z)a\Vert \le 
C\Vert a\Vert$ for all $z\in \SIGZZ$ and $a\in X$.
Therefore, by the change $s=z\eta$: $s \longrightarrow \eta$ of the 
integral variables, we obtain
$$
\Vert J^{\beta}G(z)a\Vert 
= \left\Vert \frac{1}{\Gamma(\beta)}\int^z_0 (z-s)^{\beta-1}
G(s)a ds \right\Vert 
= \left\Vert \frac{z^{\beta}}{\Gamma(\beta)}\int^1_0 (1-\eta)^{\beta-1}
G(z\eta)a d\eta \right\Vert,
$$
and so 
$$
\Vert J^{\beta}G(z)a\Vert \le C\ABSz^{\beta}
\sup_{\xi\in \SIGZZ} \Vert G(\xi)a\Vert
\le C_{\beta}\ABSz^{\beta}\Vert a\Vert
$$
for all $z\in \SIGZZ$.  Thus the proof of Lemma 4.5 is
complete.
$\blacksquare$
\\

Now we proceed to the completion of the proof of Lemma 4.2.
We have 
$$
 (-A)^{\beta}\frac{d}{dz}J^{\tau}G(z)a 
= \frac{1}{2\pi i}\int_{\Gamma} \frac{d}{dz}J^{\tau}(e^{\la z})\la^{\alpha-1}
(-A)^{\beta}(\la^{\alpha} - A)^{-1}a d\la.
                                                    \eqno{(4.5)}
$$
In this equality, the verification of the exchange of $\frac{d}{dz}J^{\tau}$ 
and $\int_{\Gamma} \cdots d\la$, is similar to (4.2), and we omit the
details.  Therefore, in terms of Lemma 3.3,  (4.3) and (4.4), we have 
\begin{align*}
& \left\Vert (-A)^{\beta} \frac{d}{dz}J^{\tau}G(z)a\right\Vert
\le \frac{1}{2\pi}\int_{\Gamma} \ABSz^{\tau-1} \vert E_{1,\tau}(\la z)\vert
\vert \la\vert^{\alpha-1}\Vert (-A)^{\beta}(\la^{\alpha} - A)^{-1}a\Vert
\vert d\la\vert                               \\
\le &C\ABSz^{\tau-1} \int_{\Gamma}\frac{1}{1+\vert \la z\vert}
\vert \la\vert^{\alpha-1}\vert \la\vert^{\alpha(\beta-1)} \vert d\la\vert
\Vert a\Vert
= C\ABSz^{\tau-1}\left( \int_{\Gamma} \frac{\vert \la\vert^{\alpha\beta-1}}
{1+\vert \la z\vert} \vert d\la\vert \right)\Vert a\Vert\\
= & C\ABSz^{\tau-1}\Vert a\Vert \left(\int_{\Gamma_1}
+ \int_{\Gamma_2} + \int_{\Gamma_3}\right) 
\frac{\vert \la\vert^{\alpha\beta-1}}
{1+\vert \la z\vert} \vert d\la\vert
=: J_1(z) + J_2(z) + J_3(z).
\end{align*}

{\bf Estimation of $J_1(t)$}.
We set $\la = \rho e^{-i\gamma}$ with $\rho > \frac{1}{\ABSz}$.
By the change $\eta:= \ABSz\rho$ and 
$\alpha\beta - 1 < 0$, we have
\begin{align*}
& \int_{\Gamma_1} \frac{\vert \la\vert^{\alpha\beta-1}}
{1+ \vert \la z\vert} \vert d\la \vert
= \int^{\infty}_{\frac{1}{\ABSz}} \frac{\rho^{\alpha\beta-1}}
{1 + \rho \ABSz} d\rho
= \int^{\infty}_1 \left( \frac{\eta}{\ABSz}\right)^{\alpha\beta-1}
\frac{1}{1+\eta}\frac{1}{\ABSz} d\eta\\
=& \ABSz^{-\alpha\beta}\int^{\infty}_1 \frac{1}{\eta^{1-\alpha\beta}}
\frac{1}{1+\eta} d\eta
\le \ABSz^{-\alpha\beta}\int^{\infty}_1 \eta^{\alpha\beta-2} d\eta
=: C\ABSz^{-\alpha\beta}.
\end{align*}
Therefore,
$$
\Vert J_1(z)\Vert \le C\ABSz^{\tau-1}\ABSz^{-\alpha\beta} \Vert a\Vert
= C\ABSz^{(\tau-\alpha\beta)-1}\Vert a\Vert.
$$

{\bf Estimation of $J_2(z)$.}
We set $\la = \frac{1}{\ABSz}e^{i\theta}$, where $\theta$ varies from 
$-\gamma$ to $\gamma$.  Then, $d\la = \frac{1}{\ABSz}ie^{i\theta}d\theta$ and 
$\vert \la\vert = \frac{1}{\ABSz}$, $\vert d\la\vert = \frac{1}{\ABSz}d\theta$.
Consequently, 
$$
\int_{\Gamma_2} \frac{\vert \la\vert^{\alpha\beta-1}}
{1+\vert \la z\vert}\vert d\la\vert
= \int_{-\gamma}^{\gamma} \left( \frac{1}{\ABSz}\right)^{\alpha\beta-1}
\frac{1}{2}\frac{1}{\ABSz} d\theta
= \frac{1}{2}\ABSz^{-\alpha\beta} (2\gamma).
$$
Hence, 
$$
\Vert J_2(z)\Vert \le C\ABSz^{\tau-1}\ABSz^{-\alpha\beta} \Vert a\Vert
= C\ABSz^{(\tau-\alpha\beta)-1}\Vert a\Vert.
$$
As for $J_3(z)$, the estimation is essentially same as $J_1(z)$.
Thus the proof of Lemma 4.2 is complete.
$\blacksquare$
\\
{\bf Third Step: Proof of (2.8) with $a=0$ and $F \in C^{\infty}_0
(0,T;\DDD(A))$.}
\\
Henceforth we write $F'(s) = \frac{dF}{ds}(s)$, $F'(\xi) = 
\frac{dF}{d\xi}(\xi)$, etc.
We will show
\\
{\bf Lemma 4.6.}
\\
{\it
Let $F \in C^{\infty}_0(0,T;\DDD(A))$.  Then
$$
\int^t_0 \left( \frac{d}{dt}J^{\alpha}G\right)(t-s)F(s) ds
\in W_{\alpha,q}(0,T;X)
$$
and
$$
\pppa \int^t_0 \left( \frac{d}{dt}J^{\alpha}G\right)(t-s)F(s) ds
= \int^t_0 G(t-s)F'(s) ds, \quad t>0.                
$$
}
\\
{\bf Proof of Lemma 4.6.}
\\
By Lemma 4.5, we see that $J^{\alpha}Ga \in L^{\infty}(0,\infty;X)$ for 
all $a \in X$.  Hence, for sufficiently small $\delta>0$, by integration 
by parts, we have
\begin{align*}
& \int^{t-\delta}_0 \left( \frac{d}{dt}J^{\alpha}G\right)(t-s)F(s) ds
= -\int^{t-\delta}_0 \frac{d}{ds}(J^{\alpha}G(t-s))F(s) ds \\
= -& [J^{\alpha}G(t-s)F(s)]^{s=t-\delta}_{s=0}
+ \int^{t-\delta}_0 J^{\alpha}G(t-s)F'(s) ds.
\end{align*}
Since $F \in C^{\infty}_0(0,T;\DDD(A))$, we have $F(0) = 0$. 
Lemma 4.5 for $z>0$ yields
$$
\Vert J^{\alpha}G(\delta)F(t-\delta)\Vert \le C\delta^{\alpha}
\Vert F(t-\delta)\Vert 
\le C\Vert F\Vert_{C([0,T];X)}\delta^{\alpha} \rrrr 0
$$
as $\delta \to 0$.  Consequently,
$$
\int^t_0 \left( \frac{d}{dt}J^{\alpha}G\right)(t-s)F(s) ds
= \int^t_0 J^{\alpha}G(t-s)F'(s) ds, \quad 0<t<T.              \eqno{(4.6)}
$$

Next, we will calculate:
$$
J^{\alpha}\left( \int^t_0 G(t-s)F'(s) ds\right)
= \frac{1}{\Gamma(\alpha)} \int^t_0 (t-s)^{\alpha-1}
\left( \int^s_0 G(s-\xi)F'(\xi) d\xi \right) ds.
$$
Exchanging the orders of the integrals:
$$
\int^t_0 \left( \int^s_0 \cdots d\xi\right) ds
= \int^t_0 \left( \int^t_\xi \cdots ds\right) d\xi,
$$
we obtain
\begin{align*}
& \frac{1}{\Gamma(\alpha)} \int^t_0 (t-s)^{\alpha-1}
\left( \int^s_0 G(s-\xi)F'(\xi) d\xi \right) ds
= \frac{1}{\Gamma(\alpha)} \int^t_0 
\left( \int^t_\xi (t-s)^{\alpha-1}G(s-\xi) ds\right) F'(\xi) d\xi\\
=& \frac{1}{\Gamma(\alpha)} \int^t_0 
\left( \int^{t-\xi}_0 (t-\xi-\eta)^{\alpha-1}G(\eta) d\eta\right) 
F'(\xi) d\xi.
\end{align*}
For the last equality, we changed the variables: $s \mapsto \eta$ by 
$\eta:= s-\xi$.

Since $\frac{1}{\Gamma(\alpha)} \int^{t-\xi}_0 
(t-\xi-\eta)^{\alpha-1}G(\eta) d\eta = (J^{\alpha}G)(t-\xi)$, 
by means of (4.6), we reach 
$$
\int^t_0 \left( \frac{d}{dt}J^{\alpha}G\right)(t-s) F(s) ds
= J^{\alpha} \left(\int^t_0 G(t-s) F'(s) ds\right), \quad 
0<t<T.                                  \eqno{(4.7)}
$$
Moreover, Lemma 3.1 (ii) with $\beta = 0$ yields
$$
\left\Vert G(t-s)F'(s) ds\right\Vert
\le \int^t_0 \Vert G(t-s)\Vert \Vert F'(s)\Vert ds 
\le C\max_{0\le s\le T} \Vert F'(s)\Vert,
$$
so that 
$$
\int^t_0 G(t-s)F'(s) ds \in L^q(0,T;X).
$$
Therefore (4.7) and the definition of $\pppa$ imply
$$
J^{\alpha}\int^t_0 G(t-s) F'(s) ds \in W_{\alpha,q}(0,T;X),
$$
that is,
$$
\int^t_0 \left( \frac{d}{dt}J^{\alpha}G\right)(t-s) F(s) ds
\in W_{\alpha,q}(0,T;X).
$$
Thus the proof of Lemma 4.6 is complete.
$\blacksquare$
\\

Next, we set 
$$
W(F)(t):= \int^t_0 \left( \frac{d}{dt}J^{\alpha}G\right)(t-s)F(s) ds
= - \int^t_0 \frac{d}{ds}(J^{\alpha}G(t-s))F(s) ds
$$
for $0<t<T$ and $F \in C^{\infty}_0(0,T;\DDD(A))$.
Then, the definition of $K$ yields
$$
W(F)(t):= \int^t_0 K(t-s)F(s) ds \quad \mbox{for $0<t<T$ and 
$F \in C^{\infty}_0(0,T;\DDD(A))$}.
                                             \eqno{(4.8)}
$$
\\
{\bf Lemma 4.7.}
\\
{\it
Let $F \in C^{\infty}_0(0,T;\DDD(A))$.  Then,
$W(F) \in W_{\alpha,q}(0,T;X)$, $AW(F) \in L^q(0,T;X)$
and
$$
\pppa W(F)(t) =  AW(F)(t) + F(t) 
\quad \mbox{for $0<t<T$}.
$$
}
\\
{\bf Proof of Lemma 4.7.}
\\
Choosing $\delta > 0$ sufficiently small, 
we set
$$
W_{\delta}(F)(t) = -\int^{t-\delta}_0 \left(\frac{d}{ds}
(J^{\alpha}G)(t-s)\right) F(s) ds.
$$
Since $F(t)\in \DDD(A)$, we obtain
$$
A\left( \frac{d}{ds} (J^{\alpha}G)(t-s)\right) F(s)
= \left(\frac{d}{ds}(J^{\alpha}G)(t-s)\right) AF(s)
\in L^1_s(0,T;X).
$$
Noting that $A$ is a closed operator in $X$, by an argument similar to 
Lemma 3.2, we see
$$
AW_{\delta}(F)(t) = -\int^{t-\delta}_0 \frac{d}{ds}
(J^{\alpha}G)(t-s) AF(s) ds
= -\int^{t-\delta}_0 \left(\frac{d}{ds}(J^{\alpha}AG)(t-s) \right) F(s) ds.
$$
Proposition 3.1 yields 
$$
AG(t-s)a = \pppa (G(t-s)a-a) \quad \mbox{in $X$ for $0<s<t$ and all 
$a\in \DDD(A)$}
$$
and
$$
J^{\alpha}AG(t-s)F(s) = J^{\alpha}\pppa (G(t-s) -1)F(s).
$$
Since $F(t)\in \DDD(A)$, we can apply (3.5) and 
$(G(t-s)-1)F(s) \in \DDD(\pppa)$ as a function in $t$.
Consequently,
$$
J^{\alpha}\pppa (G(t-s)-1) = G(t-s) -1,
$$
which implies
$$
J^{\alpha}AG(t-s)F(s) = (G(t-s)-1)F(s).
$$
Hence,
$$
AW_{\delta}(F)(t) = -\int^{t-\delta}_0 \frac{d}{ds}(G(t-s) - 1)F(s) ds,
\quad 0<t<T.
$$

Similarly to the proof of (4.6), the integration by parts implies
\begin{align*}
& AW_{\delta}(F)(t) 
= -[(G(t-s)-1)F(s)]^{s=0}_{s=t-\delta}
+ \int^{t-\delta}_0 (G(t-s)-1)F'(s) ds\\
=& (G(\delta)-1)F(t-\delta)
+ \int^{t-\delta}_0 G(t-s)F'(s) ds
- \int^{t-\delta}_0 F'(s) ds.
\end{align*}
On the other hand, as we already proved part (4) of Theorem 2.1 at the end of
Section 3, we have
$$
\lim_{t\to 0} \Vert G(t)a - a\Vert = 0            \eqno{(4.9)}
$$
for all $a \in \DDD(A)$.
Therefore, we can show  
$$
\lim_{\delta\to 0} \Vert (G(\delta)-1)F(t-\delta)\Vert 
= 0 \quad \mbox{for each $0 < t< T$.}                           \eqno{(4.10)}
$$
{\bf Proof of (4.10).}
First, 
\begin{align*}
& \Vert (G(\delta) - 1)F(t-\delta)\Vert
= \Vert (G(\delta) - 1)F(t) + (G(\delta) - 1)(F(t-\delta) - F(t))\Vert\\
\le& \Vert (G(\delta) - 1)F(t) \Vert 
+ \Vert (G(\delta) - 1)(F(t-\delta) - F(t))\Vert.
\end{align*}
Let $t>0$ be arbitrarily fixed.
Then $\lim_{\delta\to 0}\Vert (G(\delta) - 1)F(t)\Vert = 0$ by (4.9).
Moreover, Lemma 3.1 with $\beta = 0$ yields
\begin{align*}
&\Vert G(\delta)-1\Vert \Vert F(t-\delta) - F(t)\Vert 
\le C(\Vert G(\delta)\Vert + 1)\Vert F(t-\delta) - F(t)\Vert \\
\le& C\Vert F(t-\delta) - F(t)\Vert.
\end{align*}
Since $t>0$, the limit of the right-hand side is zero as 
$\delta \to 0$.
Hence $\lim_{\delta\to 0} \Vert (G(\delta) - 1)F(t-\delta)\Vert = 0$
for any fixed $t>0$.
$\blacksquare$

Therefore,
\begin{align*}
& AW(F)(t) = \lim_{\delta\to 0} AW_{\delta}(F)(t)\\
=& \lim_{\delta\to 0} (G(\delta)-1)F(t-\delta)
+ \lim_{\delta \to 0}\int^{t-\delta}_0 G(t-s)F'(s) ds
- \lim_{\delta \to 0} (F(t-\delta) - F(0)),
\end{align*}
and so
$$
AW(F)(t) = \int^t_0 G(t-s)F'(s) ds - F(t) \quad 
\mbox{for $0 < t < T$}.                     \eqno{(4.11)}
$$
On the other hand, the application of Lemma 4.6 to (4.11)
yields
$$
AW(F)(t) = \pppa W(F)(t) - F(t), \quad 0<t<T.
$$
Thus the proof of Lemma 4.7 is complete.
$\blacksquare$
\\
{\bf Fourth Step: Proof of (2.8) with $a=0$ and 
$F \in L^q(0,T; \DDD((-A)^{\ep}))$.}
\\
By (2.12), we can readily verify that 
$AW(F) \in L^q(0,T;X)$ for $F \in L^q(0,T;\DDD((-A)^{\ep}))$.
Indeed,
\begin{align*}
& \Vert AW(F)(t) \Vert
= \left\Vert \int^t_0 (-A)^{1-\ep}(-A)^{-\ep} K(t-s)F(s) ds\right\Vert\\
=& \left\Vert \int^t_0 (-A)^{1-\ep}K(t-s)(-A)^{-\ep} F(s) ds
\right\Vert
\le C\int^t_0 (t-s)^{\alpha\ep-1} \Vert \AEAE F(s)\Vert ds,\quad
0<t<T.
\end{align*}
The Young inequality on the convolution yields
\begin{align*}
& \Vert AW(F)\Vert_{L^q(0,T;X)}
= \left( \int^T_0 \Vert AW(F)(s)\Vert^q ds\right)^{\frac{1}{q}}\\
\le & C\Vert t^{\alpha\ep-1}\Vert_{L^1(0,T)}
\left( \int^T_0 \Vert \AEAE F(s) \Vert^q ds\right)^{\frac{1}{q}}
= C\Vert F\Vert_{L^q(0,T;\DDD(\AEAE))}.
\end{align*}

Next, for arbitrarily given $F \in L^q(0,T;\DDD(\AEAE))$,
by the mollifier in $t$ (e.g., Adams \cite{Ad}) and the 
density of $\DDD(A)$ in $\DDD(\AEAE)$, we can find a sequence 
$F_n \in C^{\infty}_0(0,T;\DDD(A))$, $n\in \N$ such that 
$\Vert F_n-F\Vert_{L^q(0,T;\DDD(\AEAE))} \rrrr 0$ as 
$n\to \infty$.
In terms of (4.8) and Lemma 4.7, we have
\begin{align*}
& W(F_n)(t):= \int^t_0 K(t-s)F_n(s) ds \in W_{\alpha,q}(0,T;X) \cap
L^q(0,T;X),\\
& \pppa W(F_n)(t) = AW(F_n)(t) + F_n(t), \quad 0<t<T \quad
\mbox{for $n\in \N$}.
\end{align*}

For $\Vert AW(F_n)(t) - AW(F)(t)\Vert
= \left\Vert \int^t_0 (-A)^{1-\ep}(-A)^{-\ep} K(t-s)(F-F_n)(s) ds
\right\Vert$, by means of the Young inequality on the convolution,
we can similarly estimate:
$$
\Vert AW(F_n) - AW(F)\Vert_{L^q(0,T;X)}
\le C\Vert F - F_n\Vert_{L^q(0,T;\DDD((-A)^{\ep}))}.
$$
Therefore,
$$
\Vert AW(F_n) - AW(F)\Vert_{L^q(0,T;X)} \rrrr 0
\quad \mbox{as $n\to \infty$}.                       \eqno{(4.12)}
$$
Since $\pppa W(F_n) = AW(F_n)+F_n$ in $(0,T)$, we see
$$
\pppa W(F_n) \rrrr AW(F)+F \quad \mbox{in $L^q(0,T;X)$ as 
$n\to \infty$.}                                     \eqno{(4.13)}
$$
Similarly to (4.12), we can verify 
$$
\Vert W(F_n) - W(F)\Vert_{L^q(0,T;X)} \rrrr 0 \quad 
\mbox{as $n\to \infty$}.
$$
Since $\pppa$ is a closed operator with the domain $W_{\alpha,q}(0,T;X)
\subset L^q(0,T;X)$ to $L^q(0,T;X)$ by Lemma 1.1, it follows from (4.13) that  
$W(F) \in \DDD(\pppa)$ and 
$$
\pppa W(F) = \lim_{n\to\infty} \pppa W(F_n) = AW(F)+F \quad
\mbox{in $L^q(0,T;X)$}.
$$
Therefore, (2.8) is proved for $a=0$ and $F \in L^q(0,T;\DDD(\AEAE))$.
Thus the proof of (2.8) for $F \in L^q(0,T;\DDD(\AEAE))$ is finished.
Consequently, the proof of Theorem 2.1 is complete.
$\blacksquare$
\section{Proof of Theorem 2.2 and Corollary 2.1}

{\bf 5.1. Proof of Theorem 2.2.}
\\
{\bf Proof of (i).}
By Lemma 3.4, we have
$$
L(\pppa u)(\la) = \la^{\alpha}Lu(\la) \quad \mbox{for all
$\la>0$.}
$$
Here we recall the Laplace transform $Lv(t) := \int^{\infty}_0
e^{-\la t}v(t) dt$.
Therefore, $\la^{\alpha}(Lu)(\la) = A(Lu)(\la)$ for all $\la>0$,
that is,
$$
(A-\la^{\alpha})(Lu)(\la) = 0 \quad \mbox{for all
$\la>0$.}
$$
Hence, $(Lu)(\la) = 0$ if $\la^{\alpha} \in \rho(A)$.

Since $\rho(A) \supset \Sigma_{\gamma}:= \{ z\in \C;\, 
\vert \mbox{arg}\, z\vert < \gamma,\, z\ne 0\}$ with some
$\gamma\in \left(\frac{\pi}{2}\, \pi\right)$ by 
{\it Condition ($\AAA$)}, if $\la \in \Sigma_{\gamma}$, then 
$\la^{\alpha} \in \Sigma_{\gamma} \subset \rho(A)$ by 
$0<\alpha<1$.  Therefore, $Lu(\la) = 0$ for all $\la \in 
\Sigma_{\gamma}$.
By $\{ \la > 0\} \subset \rho(A)$, we see that 
$(Lu)(\la) = 0$ for $\la > 0$.  The injectivity of the
Laplace transform yields $u(t) = 0$ for $t>0$.
Thus the proof of Theorem 2.2 (i) is completed.
$\blacksquare$
\\
{\bf Proof of (ii).}
Since $W_{\alpha,q}(0,T;X) \cap L^q(0,T;\DDD(A))
\subset W_{\alpha,2}(0,T;X) \cap L^2(0,T;\DDD(A))$, it 
suffices to prove the uniqueness for $q=2$.

First we note the coercivity (Theorem 3.4 (ii) in 
\cite{KRY}, for example):
$$
\frac{1}{\Gamma(\alpha)}\int^t_0 (t-s)^{\alpha-1}
(\ppp_s^{\alpha} u(s),\, u(s))_X ds \ge \frac{1}{2}\Vert u(t)\Vert^2_X
                            \eqno{(5.1)}
$$
for $u\in W_{\alpha,2}(0,T;X)$.

By $\pppa u = Au$ and (2.13), in terms of (5.1),
we obtain
\begin{align*}
& 0 = \frac{1}{\Gamma(\alpha)}\int^t_0 (t-s)^{\alpha-1}
(\ppp_s^{\alpha} u(s), u(s))_Xds 
- \frac{1}{\Gamma(\alpha)}\int^t_0 (t-s)^{\alpha-1}
(Au(s), u(s))_X ds\\
\ge& \frac{1}{2}\Vert u(t)\Vert_X^2 
- C_1\int^t_0 (t-s)^{\alpha-1} \Vert u(s)\Vert^2_X ds,
\quad 0<t<T.
\end{align*}
The generalized Gronwall inequality (e.g., Henry \cite{H},
Kubica, Ryszewska and Yamamoto \cite{KRY}) yields $u=0$ in $(0,T)$, 
and so the proof of  Theorem 2.2 (ii) is complete.
$\blacksquare$
\\
{\bf Proof of (iii).}
We know an estimate 
$$
\Vert v\Vert_{W^{2,p}(\OOO)} \le C(\Vert Av\Vert_{L^p(\OOO)}
+ \Vert v\Vert_{L^p(\OOO)}) \quad \mbox{for $v \in \DDD(A)$}
                      \eqno{(5.2)}
$$
(e,g., Ladyzhenskaya and Ural'tseva \cite{LU},
or Theorem 3.1 (p.212) in Pazy \cite{Pa}).
Taking into consideration that all the coefficients of 
$A$ are in $C^{\infty}(\ooo{\OOO})$, we can verify that
$A^{-\ell}a \in W^{2\ell,p}(\OOO)$ for each $\ell \in \N$ and 
$a\in L^p(\OOO)$.  Moreover, we can readily see that 
$J^mg \in W^{m,q}(0,T;L^p(\OOO))$ for each $m\in \N$ and 
$g \in L^q(0,T;L^p(\OOO))$.

In view of the Sobolev embedding,
for $1 < p < \infty$ and $1 \le q \le \infty$, we choose 
sufficiently large $\ell, m \in \N$ such that 
$$
W^{m,q}(0,T;W^{2\ell,p}(\OOO)) \subset W^{1,2}(0,T;L^2(\OOO)) \cap
L^2(0,T;H^2(\OOO)\cap H^1_0(\OOO)).
$$
Then,
$$
w:= J^mA^{-\ell}u \in 
W^{1,2}(0,T;L^2(\OOO)) \cap L^2(0,T;H^2(\OOO)\cap H^1_0(\OOO))
                                                  \eqno{(5.3)}
$$
for $u \in L^q(0,T;L^p(\OOO))$.

Noting that $J^m\pppa u = \pppa J^mu$ for 
$u \in \DDD(\pppa) = W_{\alpha,q}(0,T;L^p(\OOO))$, we operate
$J^mA^{-\ell}$ to the equation $\pppa u = Au$ in $L^q(0,T;L^p(\OOO))$,
so that we obtain
$$
\pppa w = Aw \quad \mbox{in $L^2(0,T;L^2(\OOO))$.}     \eqno{(5.4)}
$$
In terms of (5.3), for (5.4) we can apply the uniqueness of solution 
$w \in L^2(0,T;\DDD(A)) \cap W_{\alpha,2}(0,T;L^2(\OOO))$
(e.g., \cite{KRY}, \cite{Za2009}), and so $w=J^mA^{-\ell} u =0$ in 
$\OOO\times (0,T)$.
As is readily verified, the operator $J^m$ and $A^{-\ell}$ are injective, and
we reach $u=0$ in $\OOO\times (0,T)$.
Thus the proof of Theorem 2.2 (iii) is complete.
$\blacksquare$

{\bf 5.2. Proof of Corollary 2.1.}
\\
We rewrite (2.1) as
$$
\pppa (u(t)-a) = A_0u(t) + C_0u(t) + F(t), \quad 0<t<T.
$$
We can construct the operator $G(t)$ and $K(t)$ by (2.2) and (2.6) for $A_0$.  
We introduce an iteration scheme:
$$
\left\{ \begin{array}{rl}
& u_0(t) = G(t)a + \int^t_0 K(t-s)F(s) ds, \cr \\
& u_{n+1}(t) = G(t)a + \int^t_0 C_0K(t-s)u_n(s) ds
+ \int^t_0 K(t-s)F(s) ds, \quad n\in \N \cup \{0\}.
\end{array}\right.
                                  \eqno{(5.5)}
$$
By Theorem 2.1 (2), we see that $u_n$ exists and well-defined, and 
$$
u_0 \in L^q(0,T;\DDD(A_0)), \quad u_0-a \in W_{\alpha,q}(0,T;X)
                                            \eqno{(5.6)}
$$
and we can estimate
$$
\Vert A_0u_0(t)\Vert \le Ct^{-\alpha(1-\mu)}\Vert (-A_0)^{-\mu}a\Vert
+ C\int^t_0 (t-s)^{\alpha\ep-1} \Vert (-A_0)^{\ep} F(s)\Vert ds,
\quad 0<t<T.
$$
By (2.14), we have
\begin{align*}
& \Vert A_0(u_2(t) - u_1(t))\Vert
= \left\Vert \int^t_0 A_0C_0K(t-s) (u_1(s) - u_0(s))ds \right\Vert\\
= &\left\Vert \int^t_0 C_0K(t-s) A_0(u_1(s) - u_0(s))ds \right\Vert\\
\le& C\int^t_0 \Vert K(t-s)\Vert \Vert A_0(u_1(s) - u_0(s))\Vert ds 
\le C\int^t_0 (t-s)^{\alpha-1}\Vert A_0(u_1(s) - u_0(s))\Vert ds.
\end{align*}
Hence, repeatedly applying similar estimates, we have
\begin{align*}
& \Vert A_0(u_3(t) - u_2(t))\Vert
\le C\int^t_0 (t-s)^{\alpha-1}
C\left(\int^s_0 (s-\xi)^{\alpha-1} \Vert A_0(u_1(\xi) - u_0(\xi))\Vert 
d\xi \right)ds\\
=& C^2\int^t_0 \left( \int^t_{\xi} (t-s)^{\alpha-1}
 (s-\xi)^{\alpha-1} ds \right) \Vert A_0(u_1(\xi) - u_0(\xi))\Vert
d\xi\\
=& C^2\int^t_0 \frac{\Gamma(\alpha)^2}{\Gamma(2\alpha)}
 (t-s)^{2\alpha-1} \Vert A_0(u_1(\xi) - u_0(\xi))\Vert d\xi.
\end{align*}
Here we exchanged the orders of the integrals and 
$$
\int^t_{\xi} (t-s)^{\alpha-1}(s-\xi)^{\alpha-1} ds
= \int^{t-\xi}_0 (t-\xi-\eta)^{\alpha-1}\eta^{\alpha-1} d\eta,
$$
which is verified by the change of the variables $\eta = s-\xi$.

Hence, 
$$
\Vert A_0(u_3(t) - u_2(t))\Vert
\le \frac{(C\Gamma(\alpha)^2}{\Gamma(2\alpha)}
\int^t_0  (t-s)^{2\alpha-1} \Vert A_0(u_1(s) - u_0(s))\Vert ds.
$$
Continuing the estimation, we can reach 
$$
\Vert A_0(u_{n+1}(t) - u_n(t))\Vert
\le \frac{(C\Gamma(\alpha)^n}{\Gamma(n\alpha)}
\int^t_0  (t-s)^{n\alpha-1} \Vert A_0(u_1(s) - u_0(s))\Vert ds, \quad
0<t<T.
$$
Therefore, in terms of (5.6), the Young inequality on the 
convolution yields
\begin{align*}
& \Vert A_0(u_{n+1} - u_n)\Vert_{L^q(0,T;X)}
= \Vert u_{n+1} - u_n\Vert_{L^q(0,T;\DDD(A_0))} \\
\le & \frac{(C\Gamma(\alpha))^n}{\Gamma(n\alpha)}
\frac{T^{n\alpha}}{n\alpha} \Vert A_0(u_1 - u_0)\Vert_{L^q(0,T;X)}
\end{align*}
$$
\le \frac{(C\Gamma(\alpha)T^{\alpha})^n}{\Gamma(n\alpha+1)}
(\Vert A_0u_1\Vert_{L^q(0,T;X)} + \Vert A_0u_0\Vert_{L^q(0,T;X)})
                                       \eqno{(5.7)}
$$
for each $n\in \N$.
Since 
$$
\lim_{n\to \infty} 
\frac{(C\Gamma(\alpha)T^{\alpha})^{n+1}}{\Gamma((n+1)\alpha+1)}
\left( \frac{(C\Gamma(\alpha)T^{\alpha})^n}{\Gamma(n\alpha+1)}
\right)^{-1}
= C\Gamma(\alpha)T^{\alpha} \lim_{n\to \infty}
\frac{\Gamma(n\alpha+1)}{\Gamma((n+1)\alpha+1)} = 0,
$$
we see 
$$
\sum_{n=0}^{\infty} 
\frac{(C\Gamma(\alpha)T^{\alpha})^n}{\Gamma(n\alpha+1)} < \infty.
$$
It follows from (5.7) that 
$$
\left\{ \begin{array}{rl}
& \Vert A_0u_n\Vert_{L^q(0,T;X)} 
\le C(\Vert A_0u_1\Vert_{L^q(0,T;X)} + \Vert A_0u_0\Vert
_{L^q(0,T;X)}) \quad \mbox{for each $n\in \N$}, \cr\\
& A_0u_n \rrrr A_0u \quad \mbox{in $L^q(0,T;X)$ as $n\to \infty$}.
\end{array}\right.
                                    \eqno{(5.8)}
$$
By Theorem 2.1, we have
$$
\pppa (u_{n+1} - a) = A_0u_{n+1} + C_0u_n + F(t),\quad 0<t<T.
$$
In view of (2.6), we see that 
$\lim_{n\to \infty} \pppa (u_{n+1} - a)=: \www{u}$ is convergent 
in $L^q(0,T;X)$ and $\www{u} = A_0u + C_0u + F(t)$ for $0<t<T$.
By the closedness of the operator $\pppa: W_{\alpha,q}(0,T;X) 
\subset L^q(0,T;X) \rrrr L^q(0,T;X)$, we obtain
$$
\pppa (u-a) = A_0u + C_0u + F(t), \quad 0<t<T.
$$
With (5.8), the proof of Corollary 2.1 is complete.
$\blacksquare$
\\

\section{Comparison of our formulation with the formulation 
by a Volterra integral equation}

In this article, we formulate an initial value problem for time-fractional
equation as an evolution equation in a Banach space $X$:
$$
\ppp_t^{\alpha}(u(t) - a) = Au(t) + F(t), \quad 0<t<T \quad 
\mbox{in $X$}.                                   \eqno{(6.1)}
$$
For simplicity, we consider only $p=2$.
We assume that the operator $A: \DDD(A) \longrightarrow
X$ satisfies Condition ($\AAA$).
A classical formulation of the initial boundary value problem for
time-fractional diffusion equation is by means of 
a Volterra equation:
$$
u(t) = G(t) + J^{\alpha}(Au)(t), \quad 0<t<T \quad \mbox{in $X$}
                                \eqno{(6.2)}
$$
with given $G$ in some function space.

The theory based on (6.2) has been comprehensively developed 
for example by Pr\"uss \cite{Pr}.  See also the references therein
and especially Bazhlekova \cite{Ba}.
The results in Sections 1 - 3 of Chapter 1 in \cite{Pr}
hold for more general integral operators than 
$J^{\alpha}$, but we are limited to the case of $J^{\alpha}$
for comparing the formulations (6.1) and (6.2).

For all $\alpha \in (0,1)$, we can reduce our formulation (6.1) with 
$F\in L^2(0,T;X)$ directly to (6.2):
by operating $J^{\alpha}$ to (6.1) and recalling the definition and so
$\pppa (u-a) = (J^{\alpha})^{-1}(u-a)$ for $u-a \in H_{\alpha}
(0,T;X) = J^{\alpha}L^2(0,T;X)$, we can readily verify that 
if $u \in L^2(0,T;\DDD(A))$ satisfies $u-a \in H_{\alpha}(0,T;X)$ and
(6.1), then $u$ satisfies (6.2) with 
$$
G(t):= J^{\alpha}F(t) + a, \quad 0<t<T.             \eqno{(6.3)}
$$

Now we will observe some inconsistency in (6.2) against (6.1) which 
formulates an initial value $a \in X$ and a source $F \in 
L^2(0,T;X)$ for all $0<\alpha<1$.
\\
{\bf (I)}
In (6.1), considering $F\in L^2(0,T;X)$,
in terms of (6.3), based on (6.1), we can consider (6.2) for  
$$
G-a \in H_{\alpha}(0,T;X)               \eqno{(6.4)}
$$
(e.g., Gorenflo, Luchko and Yamamoto \cite{GLY},
Kubica, Ryszewska and Yamamoto \cite{KRY}, Yamamoto \cite{Y22},
\cite{Y25}).
On the other hand, Proposition 1.1 on p.33 of Pr\"uss \cite{Pr} assumes
$G \in C([0,T];X)$ for the well-posedness of the problem (6.2) in some 
solution class.
We can prove that the two sets
$C([0,T];X)$ and $H_{\alpha}(0,T;X)$ have no inclusions 
for $0<\alpha<\frac{1}{2}$, and so our framework admits a different space 
of sources $F$, that is, $L^2(0,T;X)$ (see also \cite{KRY}, \cite{Y22, Y25}).
On the other hand, if $\frac{1}{2} < \alpha < 1$, the Sobolev embedding 
yields $H_{\alpha}(0,T;X) \subset C([0,T];X)$, so that the class of
$F$ is transformed to the space of $G$ through $J^{\alpha}$.
We here emphasize that (6.1) provides a uniform formulation for 
$F \in L^2(0,T;X)$, and as seen later, (6.2) is more difficult for
$0< \alpha<\frac{1}{2}$. 
\\
{\bf (II)}
The formulation (6.2) does not well specify an initial value $a$ and 
a source term $F$ if $0 < \alpha < \frac{1}{2}$ and 
$G \not\in C([0,T];X)$, whenever we consider an initial value and 
a source for a time-fractional diffusion equation.
\\
Indeed, let (6.2) be formulated with 
$G \in H^{\alpha}(0,T;X)$ for $0<\alpha<\frac{1}{2}$.
Then, we recall that $H^{\alpha}(0,T;X) = H_{\alpha}(0,T;X)$
(e.g., \cite{KRY}, \cite{Y22, Y25}).  We assume $u\in H^{\alpha}(0,T;X)
\cap L^2(0,T;\DDD(A))$ 
satisfies (6.2).  Then, for arbitrarily chosen $a \in X$, by (6.2) we have
$u(t) - a = G(t) - a + J^{\alpha}Au(t)$, and so 
$$
u(t) - a \in H^{\alpha}(0,T;X) = H_{\alpha}(0,T;X) 
= J^{\alpha}L^2(0,T;X).
$$
Since $(J^{\alpha})^{-1} = \pppa$ in $H_{\alpha}(0,T;X)$ by our 
definition, we obtain
$$
\pppa (u-a) = Au(t) + (\pppa G(t) - \pppa a).
$$
By using  
$$
\pppa a = (J^{\alpha})^{-1}a = \frac{t^{-\alpha}}{\Gamma(1-\alpha)}a
\in L^2(0,T;X)
$$
by $J^{\alpha}\left( \frac{t^{-\alpha}}{\Gamma(1-\alpha)}a\right)
= a$, it follows that $u$ satisfies 
$$
\pppa (u(t)-a) = Au(t) + \left( \pppa G(t) 
- \frac{t^{-\alpha}}{\Gamma(1-\alpha)}a\right) \quad \mbox{in $X$ for 
$0<t<T$.}                               \eqno{(6.5)}
$$
Although we can easily see that (6.5) is equivalent to 
$$
\pppa u(t) = Au(t) + \pppa G(t) \quad \mbox{in $X$ for 
$0<t<T$,}                               \eqno{(6.6)}
$$
in (6.5) we interpret $a$ and $\left( \pppa G(t) 
- \frac{t^{-\alpha}}{\Gamma(1-\alpha)}a\right)$ as an initial value and 
a source term respectively, while in (6.6) we have to understand that 
the initial value is $0$ and the source term is $\pppa G(t)$.
In particular, in formula (2.2) in Chapter 2 of \cite{Ba}, in 
(6.2) one can interpret that $a\in X$ is an initial value:
$$
u(t) = a + J^{\alpha}Au(t), \quad 0<t<T.    \eqno{(6.7)}
$$
However, this is not a unique way, because this implies
(6.1) with $F(t):= \frac{t^{-\alpha}}{\Gamma(1-\alpha)}a$:
$$
\pppa u(t) = Au(t) + F(t), \quad 0<t<T,
$$
where the intial value is $0$ and a source term is not zero.
Indeed, $\pppa a = \frac{t^{-\alpha}}{\Gamma(1-\alpha)}a$ for 
$0<t<T$.
This means that even one solution to (6.2) may simultaneously 
admit infinitely many pairs
of initial values $a$ and sources $F$ for $0 < \alpha < \frac{1}{2}$.
Such multivocality caused by (6.2) around initial values and source terms,
is mathematically consistent.  It does not, however, enable us 
to discuss various applications such as 
\\
(i) control problems with distributed inputs in the form
of source terms $F$
\\
(ii) inverse problems such as the determination of an initial value. 

For control problems and inverse problems, our formulation (6.1) is adequate.
\\
{\bf (III)}
In \cite{Ba} and \cite{Pr}, it is mainly assumed a stronger
regularity $a \in \DDD(A)$ in (6.6) than ours, but (6.1) is well-posed 
for any $a \in X$ (e.g., Theorem 4.1 in Kubica, Ryszewska and Yamamoto 
\cite{KRY}). 
Thus, the formulation (6.2) is not feasible for discussing solutions in 
weaker classes. 
\\
{\bf (IV)}
In (II), we point out that the formulation (6.2) with 
given $G \in H^{\alpha}(0,T;X)$, does not well specify
initial values and sources for $0<\alpha<\frac{1}{2}$.   Finally, we discuss 
also the cases $\alpha= \frac{1}{2}$ and $\frac{1}{2} < \alpha < 1$.

We can prove
\\
{\bf Proposition 6.1.}
\\
{\it 
Let $\frac{1}{2} < \alpha < 1$.
\\
(i) Let $u \in L^2(0,T;\DDD(A))$ satisfy (6.2) with $G \in H^{\alpha}(0,T;X)$.
Then, $G - G(0) \in H_{\alpha}(0,T;X)$, $u - G(0)$ in $H_{\alpha}(0,T;X)$ and
$$
\pppa (u - G(0)) = Au + \pppa (G - G(0)), \quad 0<t<T.    \eqno{(6.8)}
$$
\\
(ii) If $u \in L^2(0,T;\DDD(A))$ satisfies $u-a \in H_{\alpha}(0,T;X)$ 
and (6.1) with $a \in X$ and $F\in L^2(0,T;X)$, then $u$ satisfies 
(6.2) with 
$$
G(t):= J^{\alpha}F(t) + a, \quad 0<t<T.
$$
}
\\

For $\frac{1}{2} < \alpha < 1$, Proposition 6.1 implies the equivalence of 
$$
\mbox {(6.1)} \quad \mbox{with} \quad a = G(0), \quad F=\pppa (G - G(0))
$$
and
$$
\mbox {(6.2)} \quad \mbox{with} \quad G(t) = J^{\alpha}F(t) + a
$$
within the described class of solutions.

In the rest case $\alpha = \frac{1}{2}$, it is still  
not clear how to find $a\in X$ satisfying (6.4) for given $G$, so that 
we do not know how we can interpret initial values and source terms in (6.2).
Therefore, also by the argument in (II),
we have no such equivalence between (6.1) and (6.2) for
$0 < \alpha \le \frac{1}{2}$.

To sum up, we can conclude:
In (6.2), the physical interpretation of initial values and source terms is
\begin{itemize}
\item
not possible for $0 < \alpha \le \frac{1}{2}$,
\item
smoothly done for $\frac{1}{2} < \alpha < 1$.
\end{itemize}
We further emphasize that \cite{PS} excludes the case $0<\alpha 
\le \frac{1}{2}$ (Theorem 4.5.15 on p.192) when we consider the case $p=2$,  
while our formulation (6.1) provides a unified framework for 
$0 < \alpha < 1$.
\\
\vspace{0.2cm}
\\
{\bf Proof of Proposition 6.1.}
\\ 
First, we note by $\frac{1}{2} < \alpha < 1$ that 
$H^{\alpha}(0,T;X) \subset C([0,T];X)$, and so 
$G - G(0) \in H_{\alpha}(0,T;X)$ for each $G \in H^{\alpha}(0,T;X)$.

Next (6.2) implies 
$$
u(t) - G(0) = J^{\alpha}Au(t) + (G(t) - G(0)), \quad 
0<t<T.                                               \eqno{(6.9)}
$$
Since $G - G(0) \in H_{\alpha}(0,T;X) = J^{\alpha}L^2(0,T;X)$, 
equation (6.9) yields $u - G(0) \in H_{\alpha}(0,T;X)$.
By $\pppa = (J^{\alpha})^{-1}$ on $H_{\alpha}(0,T;X)$ by the definition, 
we operate $\pppa$ to (6.9), so that we can reach (6.8). 
Thus, the proof of (i) is complete.

The proof of (ii) is straightforward by operating $J^{\alpha}$ to 
(6.1).
$\blacksquare$

\section{Conclusions}

We consider an evolution equation 
$$
\pppa (u(t) - a) = Au(t) + F(t), \quad 0<t<T             \eqno{(7.1)}
$$
in Banach space $X$ with time-fractional differential operator $\pppa$ 
of the order 
$\alpha \in (0,1)$, and $a \in X$ corresponds to an initial value.
If we choose $X=L^p(\OOO)$ with bounded domain $\OOO \subset \R^d$ and 
suitable elliptic operator $A$ attached with boundary condition, 
then (7.1) describes an initial boundary 
value problem for a time-fractional diffusion equation.
 

The main subject is the well-posedness for initial value prolems (7.1).
Our method is based on the vector-valued Laplace transform, and is 
a modification of a classical construction of the analytic semigroup.
The construction derives solution formula which is feasible and convenient.

We intend to totally transfer the analytic semigroup approach 
(i.e., $\alpha=1$) to 
time-fractional evolution equations. 
In this article, we established a theory for time-fractional evolution 
equations which is feasible for various applications.  
The applications are demonstrated in a succeeding article.
\\
\vspace{0.2cm}
\\
{\bf Funding and Conflicts of interests/competing interests}
\\
Masahiro Yamamoto was supported by 
Grant-in-Aid for Scientific Research (A) 20H00117 
and Grant-in-Aid for Challenging Research (Pioneering) 21K18142 of 
Japan Society for the Promotion of Science.
Most of this work has been carried out when Masahiro Yamamoto was 
a visiting professor at Sapienza University of Rome 
in October - November 2024.

G. Floridia was supported by the INdAM (Istituto Nazionale di Alta Matematica
"F. Severi") Group for Mathematical Analysis, Probability and Applications 
(GNAMPA). 
G. Floridia and M. Yamamoto were supported by the PRIN 2022 PNRR Project 
P20225SP98 "Some mathematical approaches to climate change and its impacts" 
(funded by the European Union - NextGenerationEU, CUP B53D23027770001), 
and by the Sapienza research project "Progetto Medio 2023" 
n. RM123188F7C2E732 "Controllability and inverse problems for PDEs with 
applications to biomathematics, climatology, and smart materials", coordinated 
by G. Floridia.
\\
{\bf Declaration} The authors 
declare no conflicts of interests and no competing 
interests.

\end{document}